\documentclass[12pt]{amsart}
\usepackage{amssymb,amsmath,amsthm}
\usepackage{epsfig}

\theoremstyle{plain}
\newtheorem{theorem}{Theorem}[section]

\theoremstyle{definition}
\newtheorem*{remark}{Remark}
\newtheorem{lemma}[theorem]{Lemma}

\title{The energy of crumpled sheets in F\"{o}ppl-von K\'arm\'an plate
theory}

\author{Shankar C. Venkataramani}

\address{Department of Mathematics,
  University of Chicago, 5734 University Ave., Chicago, IL 60637,
  USA}

\date{\today}

\begin{document}

\begin{abstract}
 We study investigate a long, thin rectangular elastic membrane that is
 bent through an angle $2 \alpha$, using the F\"{o}ppl--von K\'arm\'an
 {\em ansatz} in a geometrically linear setting. We study the
 associated variational problem, and show the existence of a minimizer
 for the elastic energy. We also prove rigorous upper and lower bounds
 for the minimum energy of this configuration in terms of the plate
 thickness $\sigma$ and the bending angle, and we also obtain results
 for the structure of the {\em elastic ridge} along it's length.
 \end{abstract}

\maketitle

\section{Introduction} \label{sec:intro}

Crumpled elastic sheets can be thought of as minimum energy
configurations for the F\"{o}ppl -- von K\'arm\'an (FvK) energy. Using
this approach of elastic energy minimization, the crumpling response
is understood as a result of the elastic energy of the sheet
concentrating on a small subset of the entire sheet
\cite{pomeau,science.paper,eric}. The energy in a crumpled sheet is
concentrated on a network of thin line-like creases (ridges) that meet
in point-like vertices. Recent work has resulted in quantitative
understanding of both the vertices
\cite{vertex,maha,dcone_exp,BPCB00,MB02} and the ridges
\cite{lobkovsky,LobWit,ridge_buckling,ridge_strength}. Scaling laws
governing the behavior of crumpled sheets have been obtained
\cite{pomeau,science.paper,lobkovsky} using scaling arguments.
  
Minimum energy configurations for the FvK energy have also been
studied in the context of the {\em blistering problem}, {\em viz.} the
buckling of membranes as a result of isotropic compression along the
boundary \cite{GiOrtz,GDeSOC02}.

There is a considerable body of mathematical work focused on the
blistering problem \cite{GiOrtz,KJ,Audoly,JS,JS1,BCDM,BCDM02}. Upper
and lower bounds have been obtained for approximations to the elastic
energy \cite{GiOrtz,KJ,JS1}, for the FvK energy \cite{JS,BCDM} and for
full three dimensional nonlinear elasticity \cite{BCDM02}. The FvK
energy and full three dimensional nonlinear elasticity give the same
scaling for the upper and the lower bounds. 

Our goal in this paper is to prove corresponding rigorous results for
the energy in a {\em minimal ridge} -- a single crease in a crumpled
sheet. In addition to scaling results for the energy, we also
investigate the structure of the ridge by obtaining pointwise bounds
for its ``width''. Our results for the ridge show an interesting
contrast with the corresponding results for the blistering problem
\cite{BCDM,BCDM02}. In particular, the scaling of the energy with the
thickness of the sheet has a {\em different} exponent. This implies
that the boundary conditions play an important role in determining the
$\Gamma$--limit of the FvK energy in the limit the thickness goes to
zero. We discuss this issue further in Sec.~\ref{sec:discussion}.

This paper is organized as follows -- In Sec.~\ref{sec:setup}, we
describe the problem of interest, set up the relevant energy
functional and determine the appropriate boundary conditions. We also
rescale the various quantities to a form that is suitable for further
analysis, and recast the problem in terms of the rescaled
quantities. In Sec.~\ref{sec:l_bnd}, we prove a lower bound for the
elastic energy for our boundary conditions. In Sec.~\ref{sec:u_bnd},
we prove the corresponding upper bounds by explicit construction of a
test solution. In Sec.~\ref{sec:structure}, we investigate the
structure of a single ridge, and we present a concluding discussion in
Sec.~\ref{sec:discussion}.

\section{The variational problem} \label{sec:setup}

\begin{figure}[htbp]
  \begin{center}
\centerline{\epsfig{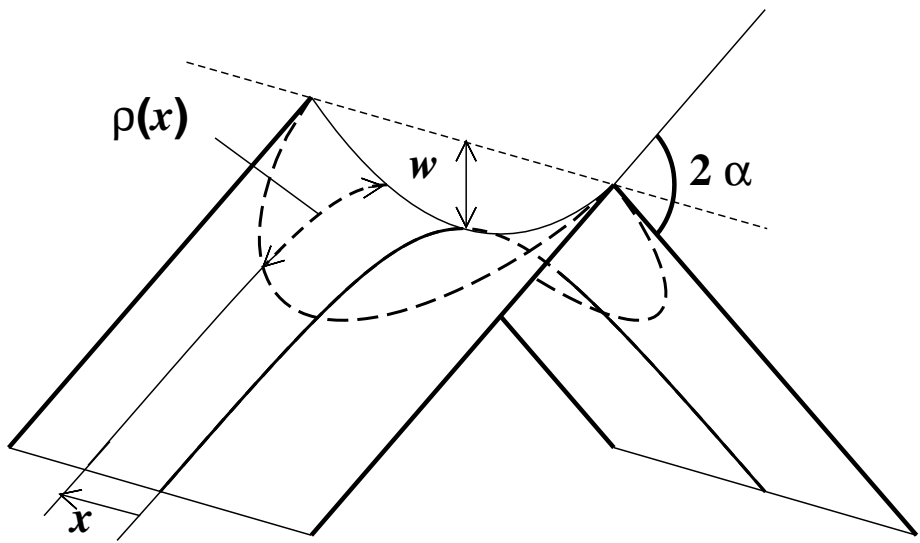}}
\end{center}
\caption{A minimal
    ridge. The boundary conditions are given by a frame (the thick
    solid lines) bent through an angle. The sheet is essentially flat
    outside the region bounded by the two dashed curves, and the bulk
    of the energy is concentrated in this region.}
\label{fig:ridge}
\end{figure}

We are interested in a {\em minimal ridge}, {\em i.e.}, the single
crease that is formed when a long rectangular elastic strip is bent
through an angle by clamping the lateral boundaries to a bent
frame. This situation is depicted in Figure~\ref{fig:ridge}.

As we will see below, the idealized boundary conditions with a sharp
corner are not appropriate, since they lead to an infinite energy for
a sheet with a finite thickness. If we make the corner extremely
sharp, all the energy (asymptotically) will be at the corners, and
this obscures the interesting physics in the problem, namely the
energy and the structure of the ridge. Thus we have to incorporate the
smoothness of the corners into our boundary conditions.  In general,
curvatures on scale smaller than the thickness cannot occur for a real
sheet, and our model energy is not appropriate for this situation. For
a crease in a real crumpled sheet, the corner is definitely smooth on
the scale of the thickness of the sheet.

We will consider the situation where the sheet is clamped to a frame,
that is much like the idealized situation depicted in
Fig.~\ref{fig:ridge}. The sheet is a rectangular strip $|x| \leq L,
|y| \leq L'$. We generally consider the situation $L \lesssim L'$. The
sheet is clamped to the frame at $x = \pm L$ and at $y = \pm L'$. We
will place the following requirements on the frame --
\begin{enumerate}
\item The boundary conditions at the frame are {\em non-stretching},
{\em i.e}, the strain $\gamma_{yy}$ is identically zero on the
boundaries $x = \pm L$, and the strains $\gamma_{xx} = \gamma_{xy} =0$
at $y = \pm L'$. 
\item The bending at the boundaries is {\em localized}. If $s$ is an
arclength parameter along the boundary, we will require that both the
boundaries are straight on a set of the form $|s| \geq k$, where $k
\ll L$. 
\item We will assume that the two ``bent'' boundaries on the frame ($x
= \pm L$) are given by {\em planar curves}, and that the angle between
the straight sections on these boundaries for $y > k$ and $y < -k$ is
$2 \alpha$ for {\em both} boundaries.
\item We will require that the straight segments for both the ``bent''
boundaries are parallel and a distance $L$ apart. This implies that
the two boundaries are not {\em twisted} with respect to one another.
\item We will also require that the two boundaries have the same {\em
asymptotic shift}. This is necessary to make $\gamma_{xy} = 0$ at $y =
\pm L'$. To define this precisely, we need to introduce appropriate
coordinates, and we will do this below.
\end{enumerate}

\subsection{Coordinates} 

We will use (the material) coordinates $(x,y)$ on the reference half
strip $|x| \leq L$, $|y| \leq L'$.  The planes containing the parallel
straight portions of the two boundaries pick out preferred {\em
in-plane} directions for $y \gg k$ and $y \ll -k$. We introduce two
sets of {\em Cartesian} coordinate systems in the ambient space. The
coordinate direction $u$ is perpendicular to the straight portion of
the boundary at $x = -L$ , and is directed toward the boundary at $x =
L$. The coordinate directions $v^+$ and $v^-$ are along the straight
portions of the boundaries for $s > k$ and $s < -k$ respectively. The
coordinate directions $w^+$ and $w^-$ give the {\em out of plane}
directions, and are chosen so that $(u,v^+,w^+)$ and $(u,v^-,w^-)$ are
right handed orthogonal triads. Finally, the origins of both the
coordinate systems coincide, and they are chosen such that the
straight portions of the boundaries lie in the plane $w^{\pm} = 0$,
and the frame boundaries are in the planes $u = \pm L$. Note that this
completely specifies the definition of the coordinate systems, and in
particular, we do not have a freedom to translate $v^{\pm}$.

\begin{figure}[htbp]
  \begin{center}
\centerline{\epsfig{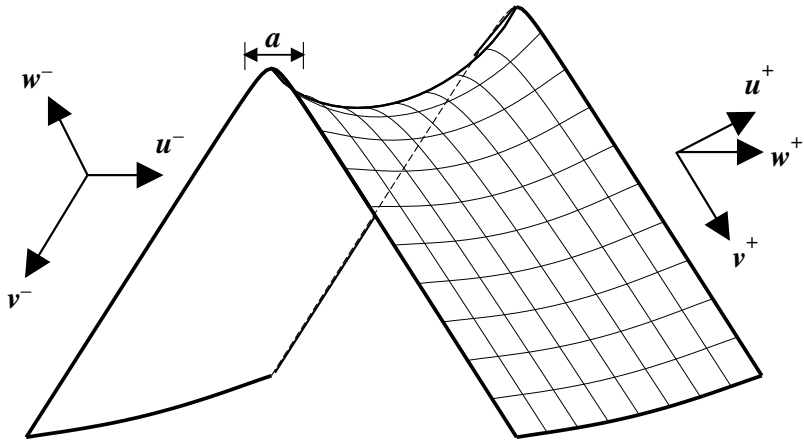}}
\end{center}
\vspace{0cm}
\caption{A schematic representation of our coordinate system and the
boundary conditions imposed on the sheet. The thick lines depict the
``frame'', and the ``corners'' are smooth on a scale $a$. $u^{\pm}$
and $v^{\pm}$ are the in-plane directions and $w^{\pm}$ are the
out-of-plane directions. The grid is given by the lines $x =
\mbox{constant}$ and $y = \mbox{constant}$}.
\label{fig:coords}
\end{figure}

The various coordinate are represented schematically in
Figure~\ref{fig:coords}. The grid in the figure is generated by the
lines $x = \mbox{constant}$ and $y = \mbox{constant}$ that are
straight in the reference (material) coordinates. We will use the
coordinate system $(u,v^+,w^+)$ for the portion of the sheet with $y
\geq 0$ and $(u,v^-,w^-)$ for $y \leq 0$. At $y = 0$, we have the
matching conditions
\begin{equation} \left(\begin{array}{c} w^- \\ v^-\end{array} \right) =
\left(\begin{array}{cc} \cos 2 \alpha & - \sin 2 \alpha \\ \sin 2
\alpha & \cos 2 \alpha \end{array} \right) \left(\begin{array}{c} w^+
\\ v^+\end{array} \right)
\label{eq:pos_match}
\end{equation}
Also, in the straight portion of the boundaries $|s| \geq k$, since
$\gamma_{yy} = 0$, it follows that $v^{\pm}(\pm L,y) - y$ is a
constant for sufficiently large (small) $y$. For the boundary at $x =
-L$, we define the {\em asymptotic shifts} by
$$
\delta^{\pm}_1 = v^{\pm}(- L,y) - y, \quad \quad y \geq k \quad (\mbox{respectively } y \leq -k).
$$ Similarly, the asymptotic shifts for the boundary at $x = L$ are
constant if $|y|$ is sufficiently large. For the boundary at $x =
-L$, we define the asymptotic shifts by
$$ \delta^{\pm}_2 = v^{\pm}(L,y) - y, \quad \quad y \geq k \quad
(\mbox{respectively } y \leq -k).
$$ We will say that the two boundaries are {\em compatible}, if
$\delta^+_{1} = \delta_2^+$ and $\delta^-_1 = \delta^-_2$. This
clearly a necessary condition for the existence of a configuration of
the sheet that satisfies the boundary conditions, and is {\em
asymptotically strain free}, {\em i.e}, the strain is identically zero
for $|y| \geq l$ for a sufficiently large $l$. In particular, we can
take $l = k$.

Assuming that the two boundaries are compatible, we will set $\delta^+
= \delta^+_{1} = \delta_2^+$, $\delta^- = \delta^-_1 =
\delta^-_2$, and $\delta = \delta^{+} - \delta^{-}$. The quantities
$\delta^+$ and $\delta^-$ change under translations of the coordinate
$y$, but $\delta$ is an invariant under these translations, and is
purely a geometrical property of the frame. Below, we will give an
expression for $\delta$ in terms of functions specifying the boundary
conditions. W.L.O.G, we can, and henceforth will, translate the
coordinate $y$ such that $\delta^+ = - \delta^- = \frac{1}{2} \delta$.

\subsection{The Elastic energy}

A mathematically justified way to obtain the elastic energy of the
deformed sheet is to treat the sheet as a three dimensional (albeit
thin) object and use a full nonlinear three dimensional elastic energy
functional for the energy density. The sheet is now a three
dimensional object $\mathcal{S} \times
\left[-\frac{h}{2},\frac{h}{2}\right]$ with thickness $h$. Let $x,y$
denote the in plane coordinates as above, and $z$ denote the
coordinate in the thin direction. If the configuration of the sheet is
given by a mapping $\phi: \mathcal{S} \times
\left[-\frac{h}{2},\frac{h}{2}\right] \rightarrow \mathbb{R}^3$ the
elastic energy is given by
$$
\mathcal{I}_{3D} = \iint_{\mathcal{S}}dx dy \int_{-h/2}^{h/2} dz W_{3D} (\nabla \phi)
$$ 

This approach however does not take advantage of the ``thinness'' of
the sheet. In particular, we would like to treat the thin sheet as a
two dimensional object. The derivation of reduced dimensional
descriptions of thin sheets has a long history. There is a classical
theory for thin elastic sheets built on the work of Euler, Cauchy,
Kirchoff, F\"{o}ppl and Von K\'arm\'an
\cite{love,landau_elastic,ciarlet2}. 


In the classical F\"{o}ppl -- von K\'arm\'an {\em ansatz}, the behavior of
the deformation $\phi$ is completely determined by the behavior of the
center-plane $z = 0$. An asymptotic expansion with this ansatz
\cite{ciarlet2} yields an effective 2-D elastic energy
$$ 
\mathcal{I} = h \left[\iint_{\mathcal{S}}dx dy W_{2D} (\phi_x, \phi_y) + 
q h^2 \iint_{\mathcal{S}} dx dy |\nabla \nu|^2\right],
$$ 
where $q$ is a nondimensional $O(1)$ factor, $\nu = \phi_x \times
\phi_y /\|\phi_x \times \phi_y \|$ is the normal to the center
surface, and $W_{2D}$ is an effective two dimensional energy. This is
the {\em geometrically nonlinear} F\"{o}ppl-- von K\'arm\'an energy of the
thin sheet.

The functional $W_{2D}$ is zero if $(\phi_x,\phi_y) \in O(2,3)$, the
set of matrices that give isometric linear mappings $\mathbb{R}^2
\rightarrow \mathbb{R}^3$. We will demand that 
$$
W_{2D} (\phi_x, \phi_y) \geq c \,\mathrm{dist}^2((\phi_x,\phi_y), O(2,3)).
$$ and that $W_{2D} (\phi_x, \phi_y) \leq \|\phi_x\|^4 + \|\phi_y\|^4$
for large $\phi_x,\phi_y$. These conditions are identical to the
conditions on the energy in Ref.~\cite{BCDM02}. 

A typical (or canonical) energy functional which satisfies these
conditions, and has the natural invariances for the problem, {\em
viz.}, the action of $O(3)$ on $(u,v,w)$ and $O(2)$ on $(x,y)$ is
\begin{align*} 
W_{2D} & = \mathrm{dist}^2((\phi_x,\phi_y), O(2,3)) \\
& = (u_x^2 + v_x^2
+ w_x^2 - 1)^2 + 2 (u_x v_x + u_y v_y + w_x w_y)^2 + (u_y^2 + v_y^2 +
w_y^2 - 1)^2
\end{align*} 
For the most part, we will work with the energy that we get by
linearizing the above expression, only in the {\em in plane}
deformation $(u,v)$ about the reference state $u_x = 1, v_y = 1$ and
all the other derivatives are zero. This yields the linearized energy
$$ 
\hat{W}_{2D} = \left[ 2(u_x - 1) + w_x^2 \right]^2 + \left[ 2(v_y -
1) + w_y^2 \right]^2 + 2 \left[ u_x + v_y + w_x w_y \right]^2.
$$ 
Note that, by linearizing the energy functional, we have destroyed
the natural invariances of the energy, and have picked out preferred
{\em in--plane} and {\em out of plane} directions in the ambient
space \cite{ciarlet2}. We will also linearize the normal vector so that
$$
\nu \approx -w_x \hat{\mathbf{e}}_u -w_y 
\hat{\mathbf{e}}_v + \hat{\mathbf{e}}_w 
$$ Using this expression for the normal vector, gives the {\em
geometrically linear} F\"{o}ppl -- von K\'arm\'an energy
\begin{align*}
\mathcal{I}_{lin} = & h \left[\iint_{\mathcal{S}} \left[ 2(u_x - 1) +
w_x^2 \right]^2 + \left[ 2(v_y - 1) + w_y^2 \right]^2 \right. \\
& \left.+ 2 \left[ u_x +
v_y + w_x w_y \right]^2 + q h^2 \left[w_{xx}^2 + 2 w_{xy}^2 + w_{yy}^2
\right]\right]
\end{align*} 
We will henceforth normalize the energy by $4h$, and also define the
scaled thickness $\sigma = \sqrt{\frac{q}{2}} h$. The linearized FvK
energy for our problem is ${\mathcal I} = {\mathcal I}^+ + {\mathcal
I}^-$, where
\begin{eqnarray}
 {\mathcal I}^+ & = & \int_{y = 0}^{L'}\int_{x =
-L}^{L}\left[\left(u_x + \frac{1}{2}{w}_x^2 - 1\right)^2 + \frac{1}{2} ({v}_x + u_y + 
{w}_x {w}_y)^2 \right.\nonumber \\
& & \left. + \left({v}_y + \frac{1}{2}{w}_y^2 - 1\right)^2  + \frac{\sigma^2}{2} ({w}_{xx}^2 + 2 {w}_{xy}^2 +
{w}_{yy}^2 ) \right]\, dx dy,
\label{eq:unscaled} 
\end{eqnarray}
where we have suppressed the subscript $+$ on $v^+$ and
$w^+$. ${\mathcal I}^-$ is given by the corresponding expression for
$y \leq 0$.


As defined above $x$ and $y$ are reference coordinates on the sheet,
$u$ and $v$ are in-plane coordinates, $w$ is the out of plane
displacement and $\sigma$ is the scaled thickness of the sheet. The
integrand includes the squares of the linearized strains,
$$ 
\gamma_{xx} = u_x + \frac{1}{2}w_x^2 - 1, \quad \gamma_{xy} =
\gamma_{yx} = \frac{1}{2} (v_x + u_y + w_x w_y), \quad \gamma_{yy} =
v_y + \frac{1}{2}w_y^2 - 1.
$$ 

\subsection{Boundary and matching conditions}

The blistering of thin films is also described by the elastic
energy in (\ref{eq:unscaled}). A similar energy also describes
multiple scale buckling in {\em free} elastic sheets ({\em i.e.}
sheets that are not forced through the boundary conditions) that are
not intrinsically flat \cite{eran,eran2}.

The difference between the blistering problem and a minimal ridge is
in the boundary conditions, which we describe below. If the bending
half-angle $\alpha \ll 1$, $\tan \alpha \approx \sin \alpha \approx
\alpha$. In this case, the deformations and the linearized strains are
small.

If $s$ is an arclength parameter along the boundary, with our choice
of coordinate systems, the frame is given by functions $(v^{\pm}_1(s),
w^{\pm}_1(s))$ for $s \geq 0$ (resp. $s \leq 0$) at $x = -L$, and
functions $(v^{\pm}_2(s), w^{\pm}_2(s))$ for $s \geq 0$ (resp. $s \leq
0$) at $x = L$. We define the {\em curvature} $\kappa_1$ of the
boundary at $x = -L$ by
$$
\kappa_1 = \frac{1}{a_1} = \int_0^{L'} ({w^+}_{ss})^2 ds + \int^0_{-L'} 
({w^-}_{ss})^2 ds
$$ $a_1$ is a length scale associated with the boundary at $x =
-L$. We define $k_1$ by demanding that the boundary at $x = -L$ is
straight for $|s| \geq k_1$. We will require that $k_1$ be on the same
scale as $a_1$, so that $k_1 = K a_1$, with $K$ staying order 1 as we
change the parameters in the problem. We can similarly define $a_2$
and $\kappa_2$ using the boundary at $x = L$. We will assume that the
length scales $a_1$ and $a_2$ are comparable, {\em i.e} $a_2/a_1 \sim
O(1)$ and likewise for $k_1$ and $k_2$. In fact, we will typically
assume that they are equal $a_1 = a_2 = a$, and $k_1 = k_2 = Ka$.

This length scale $a$ sets the natural length scale for the
boundaries. We will assume that the boundary conditions are
compatible, and further, the asymptotic shift is on the scale
$a\alpha^2$, {\em i.e.}, $\delta = \Delta a \alpha^2 $, where $\Delta$
stays order 1 as we change the parameters of the problem.

Since the linearized strain $\gamma_{yy}$ is zero at the boundaries,
we have $v^{\pm}_y = 1 - \frac{1}{2}(w^{\pm}_{y})^2$. Outside $|y| \leq k$,
$w^{\pm}_y = 0$. Since the asymptotic shift is $\delta$, we get
$v^{\pm}(\pm L,y) - y = \pm \frac{1}{2} \delta$ for $|y| > k$. 

The elastic energy penalizes the square of the
curvature. Consequently, all finite energy configurations of the sheet
have a continuous tangent plane a.e. This yields the matching
condition
$$ \left(\begin{array}{c} \nabla w^- \\ \nabla v^-\end{array} \right) =
\left(\begin{array}{cc} \cos 2 \alpha & - \sin 2 \alpha \\ \sin 2
\alpha & \cos 2 \alpha \end{array} \right) \left(\begin{array}{c} \nabla w^+
\\ \nabla v^+\end{array} \right)
$$ 
at $y = 0$, where $\nabla$ denotes the 2 dimensional
gradient. Since (\ref{eq:pos_match}) holds for all $x$ at $y = 0$, it
automatically implies that
$$ \left(\begin{array}{c} w^-_x \\ v^-_x \end{array} \right) =
\left(\begin{array}{cc} \cos 2 \alpha & - \sin 2 \alpha \\ \sin 2
\alpha & \cos 2 \alpha \end{array} \right) \left(\begin{array}{c} w^+_x
\\ v^+_x\end{array} \right)
$$ 
We therefore have an independent matching condition
\begin{equation}
\left(\begin{array}{c} w^-_y \\ v^-_y \end{array} \right) =
\left(\begin{array}{cc} \cos 2 \alpha & - \sin 2 \alpha \\ \sin 2
\alpha & \cos 2 \alpha \end{array} \right) \left(\begin{array}{c} w^+_y
\\ v^+_y\end{array} \right) \quad \quad y = 0
\label{eq:deriv_match}
\end{equation}

For small $\alpha$, we can linearize the matching conditions, and to
first order in $\alpha$, we obtain 
\begin{eqnarray}
w^-(x,0) & = & w^+(x,0) \nonumber \\
w_y^-(x,0) & = & w_y^+(x,0) - 2 \alpha  \nonumber \\
v^-(x,0) & = & v^+(x,0) + 2 \alpha  w^+(x,0) 
\label{eq:match} 
\end{eqnarray}
If we assume that $w^{\pm}_y(x,0) \sim O(\alpha)$, then for small
strains, $v_y^{\pm}(x,0) \approx 1$, and the matching condition for
$v_y$ is satisfied as an identity to first order in $\alpha$.

Since $a$ is the natural length scale at the boundaries, $w^{\pm}_y
\sim O(\alpha)$ near $y = 0$,$v^{\pm} = y \pm \frac{1}{2} \delta$ and
$v^{\pm}_y + \frac{1}{2}(w^{\pm}_y)^2 - 1 = 0$, we define the boundary conditions
in terms of {\em scaling functions} --
$$ w_{1,2}^{\pm}(s) = \alpha a
\sqrt{2} \Phi_{1,2}^{\pm}\left(\frac{s}{a}\right), \quad \quad v_{1,2}^{\pm}(s)
= \alpha^2 a \Psi_{1,2}^{\pm}\left(\frac{s}{a}\right) + s \pm a \alpha^2
\frac{\Delta}{2}.
$$ The {\em scaling functions} $(\Phi_{1,2}^{\pm}(\eta),
\Psi_{1,2}^{\pm}(\eta))$ are normalized (have curvature 1)
$$
\int_0^{L'/a} [{\Phi^+_i}'']^2 d \eta + \int^0_{-L'/a} 
[{\Phi^-_i}'']^2 d\eta  =  1, \quad \quad \mbox{for } i = 1,2,
$$ 
are strain free, 
$$
[{\Phi_i^{\pm}}']^2 + {\Psi_i^{\pm}}' = 0,  \quad \quad \mbox{for } i = 1,2,
$$
and satisfy the matching conditions 
\begin{eqnarray}
\Phi_{i}^+(0) & = & \Phi_{i}^-(0) \quad \quad \mbox{for } i = 1,2,
\nonumber \\ {\Phi_{i}^-}'(0) & = & {\Phi_{i}^+}'(0) - \sqrt{2} \quad \quad
\mbox{for } i = 1,2, \nonumber \\ {\Psi_{i}^-}(0) & = &
{\Psi_{i}^+}(0) + 2 \sqrt{2} \Phi_{i}^+(0) + \Delta \quad
\quad \mbox{for } i = 1,2, \nonumber
\end{eqnarray}

Finally, we also have the conditions $v = \pm \frac{1}{2} \delta, u =
w = 0$ at $y = \pm L'$, where $\delta$ is the asymptotic shift
determined by the (compatible) frame boundaries at $x = \pm L$.  We
will explore the relation between the frame $(\Phi^{\pm},\Psi^{\pm})$
and the value of $\Delta$ (and consequently the value $\delta$)
now. Using $w^{\pm} = 0$ for $|y| \geq k$, and the definition of
$\Delta$, we see that
$$
\Phi_{i}^{\pm}(\eta) = \Psi_{i}^{\pm}(\eta) = 0  \quad
\quad \mbox{for } |\eta| \geq K = k/a, \, i = 1,2.
$$ Using $\gamma_{yy} = 0$, we get
$$
\Psi^{\pm}_i(\eta) = -\int_{\pm L'/a}^{\eta} [{\Phi_i^{\pm}}']^2 d \eta \quad \quad \mbox{for } \eta \geq 0 (\mbox{ resp. } \eta \leq 0).
$$ 
Using this relation in the matching condition, we obtain
$$ 
\Delta = - \left[2 \sqrt{2} \Phi_{i}^+(0) + \int_{- L'/a}^0
[{\Phi_i^{-}}']^2 d \eta + \int^{L'/a}_0 [{\Phi_i^{+}}']^2 d
\eta\right].
$$ 
This gives the explicit relation between $\Delta$ and the out of
plane displacement of a no-stretching profile.

\subsection{Rescalings} 

From Lobkovsky's results \cite{lobkovsky} and the analysis in
Re.~\cite{ridge1}, we know that the dominant energy balance is between
the longitudinal curvature $w_{yy}^2$ and the transverse strain
$w_x^4$. The natural length scale in $X$ is $L$, and as we saw above,
the natural scale for $w_y$ is $\alpha$, the bending angle.

Let $L_q$ denote the natural scale for the quantity $q$. For example,
we have $L_x = L$. The above considerations lead to the conclusions
$$ 
w_y \sim \frac{L_w}{L_y} \sim \alpha \quad \quad \sigma^2 w_{yy}^2 \sim
\frac{\sigma^2 L_w^2}{L_y^4} \sim w_x^4 \sim \frac{L_w^4}{L_x^4}.
$$
This motivates the introduction of rescaled coordinates and displacements by
$$
x = L X, \quad \quad y = \sigma^{1/3} L^{2/3} \alpha^{-1/3} Y \equiv L_y Y,
$$
and 
\begin{eqnarray}
w^{\pm} & = & \sqrt{2} \sigma^{1/3} L^{2/3} \alpha^{2/3} W^{\pm} \equiv L_w W^{\pm}, \nonumber \\
v^{\pm} & = & y \pm \frac{\delta}{2} + \sigma^{1/3}
L^{2/3} \alpha^{5/3} V^{\pm} \nonumber \\
u & = & x + \sigma^{2/3} L^{1/3} \alpha^{4/3} U. \nonumber
\end{eqnarray} 
Here $\delta$ is the asymptotic shift, and is determined by the
frame boundaries. Since $\sigma, L,x, y, u, v,w, \delta$ all have
dimensions of a length, and $\alpha$ is dimensionless, it is clear
that the rescaled quantities $X,Y,U,V^{\pm},W^{\pm}$ are all
dimensionless. Note that these rescalings are {\em different} from the
rescalings in Ref.~\cite{ridge1}.

With these rescalings, the dimensionless energy $I^{\pm} = \sigma^{-5/3}
L^{-1/3} \alpha^{-7/3} {\mathcal I}^{\pm}$ is given by
\begin{eqnarray}
I^{\pm}(U,V,W) & = & \iint \left[(U_X + W_X^2)^2 +
\frac{1}{2}\epsilon^{-2/3}(V_X + U_Y + 2 W_X W_Y)^2 \right. \nonumber
\\ 
& & +
\left. \epsilon^{-4/3}(V_Y + W_Y^2)^2 \right] + 
\left[ W_{YY}^2 + 2 \epsilon^{2/3} W_{XY}^2 +\epsilon^{4/3}
W_{XX}^2 \right]dX dY
\label{eq:scaled} 
\end{eqnarray}
where we have to use the appropriate $V^{\pm}$ and $W^{\pm}$. We have
suppressed the superscripts $\pm$ on $V$ and $W$ for
clarity. $\epsilon = \sigma/(L \alpha)$ is the natural {\em small
dimensionless parameter} in the problem. We think of $\epsilon$ as a
dimensionless thickness.

Our quest for rigorous scaling results for the energy ${\mathcal I}$
reduces to the following -- Show that the rescaled energy $I = I^+ +
I^-$, of a minimizer $(U^*,V^*,W^*)$, is bounded above and below by
positive constants {\em uniform} in the dimensionless thickness
parameter $\epsilon$, as $\epsilon \rightarrow 0$.

With this rescaling, we have the matching conditions
\begin{eqnarray}
W^+(X,0) & = & W^-(X,0) \nonumber \\
W_Y^+(X,0) & = & W_Y^-(X,0) + \sqrt{2}  \nonumber \\
V^-(X,0) & = & V^+(X,0) + 2 \sqrt{2} W^+(X,0) + \frac{\delta}{L_v} \nonumber \\
& = & V^+(X,0) + 2 \sqrt{2} W^+(X,0) + A \Delta 
\label{eq:lin_match} 
\end{eqnarray}
and the boundary conditions
\begin{eqnarray}
W^{\pm}(-1,Y) = A \Phi_{1}^{\pm}\left(\frac{Y}{A}\right), & \quad   U(-1,Y) = 0, \quad &
V^{\pm}(-1,Y)
= A \Psi_{1}^{\pm}\left(\frac{Y}{A}\right), \nonumber \\
W^{\pm}(1,Y) = A \Phi_{2}^{\pm}\left(\frac{Y}{A}\right), & \quad   U(1,Y) = 0, \quad &
V^{\pm}(1,Y)
= A \Psi_{2}^{\pm}\left(\frac{Y}{A}\right), \nonumber \\
W^{\pm}(X,\pm L'/L_y) = 0, & \quad   U(X,\pm L'/L_y) = 0, \quad &
V^{\pm}(X,\pm L'/L_y) = 0,
\label{eq:bc} 
\end{eqnarray} 
where $A = a/L_y$ is the nondimensionalized inverse curvature, and
$\Phi^{\pm}_i,\Psi_i^{\pm}$ are the scaling functions from above.

\section{Existence of a minimizer} \label{sec:exist} 

From this point forward, we will use $c,C,C',C_1,C_2,$ {\em etc.} to
denote constants whose precise numerical values are not
important. These constants can change from one line to the next. By
doing this, we can focus on the scalings of the various quantities,
without worry about the numerical values of the constants in the
scaling relations. If these constants depend on a parameter $q$, we
show this dependence by writing $C(q)$ or $C_q$. Also, we will suppress
the superscripts $\pm$ whenever this does not cause any confusion.

\begin{lemma}
 $u,v^{\pm} : \mathcal{S}^{\pm} = [-L,L] \times [0,\pm L'] \rightarrow
 \mathbb{R}$ are $H^1$ functions such that $u = 0$ at $ x =
 \pm L$ and $u = 0$ at $y = \pm L'$. 
The functional
$$
J(u,v^+,v^-) = \iint_{\mathcal{S}^+} u_x^2 + (v^+_y)^2 + \frac{1}{2}(v^+_x + u_y)^2
+ \iint_{\mathcal{S}^-} u_x^2 + (v_y^-)^2 + \frac{1}{2}(v^-_x + u_y)^2
$$ 
is {\em coercive}, in the sense $\exists c,C > 0$ such that
$$ 
J(u,v^+,v^-) \geq c(\|D u\|_{L^2}^2 + \|D v^+\|^2_{L^2}+ \|D
v^-\|^2_{L^2}) - C \left[ \int [v^+(x,0) - v^-(x,0)]^2 dx\right]
$$
\end{lemma}

\begin{proof} 

We first consider the integral over $\mathcal{S}^+$. Using the
boundary conditions, and repeated integration by parts yields,
\begin{eqnarray}
\iint u_y(x,y) v^+_x(x,y) dx dy & =  & 
- \int_{-L}^{L} u(x,0) v^+_x(x,0) dx  \nonumber \\
& & - \iint u(x,y) v^+_{xy}(x,y) dx dy \nonumber \\ 
& =  &  -\int_{-L}^{L} u(x,0) v^+_x(x,0) dx \nonumber \\ 
& & +\iint u_x(x,y) v^+_{y}(x,y) dx dy 
\end{eqnarray}
Adding a similar result for $\mathcal{S}^-$, we obtain
\begin{eqnarray}
J & = & \frac{1}{2} (\|D u\|_{L^2}^2 + \|D v^+\|^2_{L^2}+ \|D
v^-\|^2_{L^2}) 
+ \int_{-L}^{L} u(x,0) \left[v^-_x(x,0)
- v^+_x(x,0)\right] dx \nonumber \\ & & +
\frac{1}{2}\iint_{\mathcal{S}^+} (u_x(x,y) + v^+_{y}(x,y))^2 dx dy +
\frac{1}{2}\iint_{\mathcal{S}^-} (u_x(x,y) + v^-_{y}(x,y))^2 dx dy
\nonumber \\ & \geq & \frac{1}{2} (\|D u\|_{L^2}^2 + \|D
v^+\|^2_{L^2}+ \|D v^-\|^2_{L^2}) \nonumber \\ & & 
- \int \left[\theta u^2(x,0) + C_{\theta} \left|v^-_x(x,0) -
v^+_x(x,0)\right|^2 \right]dx
\label{eq:bound}
\end{eqnarray}
where $\theta > 0$ can be arbitrarily small. By the trace theorem
\cite{ziemer}, and the boundary conditions $u = 0$ at $x = \pm L$ and
$y = \pm L'$, we have
$$
\int_{-L}^L u^2(x,0) dx \leq C\|Du\|_{L^2}^2.
$$ Using this inequality in the estimate~(\ref{eq:bound}) for a
sufficiently small $\theta$ we get 
$$ 
J(u,v^+,v^-) \geq c(\|D u\|_{L^2}^2 + \|D v^+\|^2_{L^2}+ \|D
v^-\|^2_{L^2}) - C \alpha^2 \left[ \int [v^+_x(x,0)-v^-_x(x,0)]^2 dx\right],
$$
and we can take $c = 1/4$ if we so choose.
\end{proof}

We are now in a position to prove the existence of a minimizer for the
function $\mathcal{I}$ with $\sigma > 0$, for our ``no
stretch'' boundary conditions.

\begin{theorem} If $\sigma > 0$, every minimizing sequence
$(u_j,v_j^{\pm},w_j^{\pm})$ for the energy functional
$\mathcal{I}$, that satisfies the boundary and the matching
conditions, has a subsequence that converges in $H^1(\mathcal{S})
\times H^1(\mathcal{S}) \times H^2(\mathcal{S})$ to a global minimizer
of $\mathcal{I}$.
\end{theorem}

\begin{proof} This conclusion follows easily from the direct method in the 
Calculus of variations \cite{evans2,Young,dacorogna}.

$(u_j,v^{\pm}_j,w^{\pm}_j)$ is a minimizing sequence for
$\mathcal{I}$ that satisfies all the boundary conditions and
the matching conditions at $y = 0$. For all $j$, we have the boundary
conditions $w^{\pm}_j = \sqrt{2}\alpha a \Phi^{\pm}_{1,2}(y/a)$ at $x = \pm L$
and $w^{\pm}_j = 0$ at $y = \pm L'$. These boundary conditions, along
with $\mathcal{I} \geq \sigma^2 \|D^2 w\|_{L^2}^2$, imply
that $w_j$ is a bounded sequence in $H^2$.

Consequently, up to extraction of a subsequence, $w^{\pm}_j
\rightharpoonup {w^{\pm}}^*$ in $H^2$. The compactness of the
embedding $H^2 \rightarrow W^{1,4}$ implies that $(w^{\pm}_x)^2,
(w^{\pm}_y)^2$ converge strongly to $({w^{\pm}_x}^*)^2$ and
$({w^{\pm}_y}^*)^2$ respectively in $L^2$. 

The standard trace theorems \cite{ziemer} imply that 
$$
{w^{\pm}}^*(\pm L,y) = \sqrt{2}\alpha a \Phi^{\pm}_{1,2}\left(\frac{y}{a}\right), \quad \quad 
{w^{\pm}}^*(x,\pm L')  =  0.
$$
Consequently, $(u_j,v^{\pm}_j,{w^{\pm}}^*)$ satisfy the no-stretch
boundary conditions for all $j$. In general, they {\em do not} satisfy
the matching conditions at $y = 0$. Rather, we have the relation
$$
v^+_j(x,0) - v^-_j(x,0) = 2 \alpha w^+_j(x,0).
$$
We set $ \tilde{u}_j = u_j -x$.
We will henceforth drop the subscript $j$ where it wont cause any
confusion. Since $(u_j,v^{\pm}_j,{w^{\pm}}^*)$ satisfy the no-stretch
boundary conditions, $\tilde{u} = 0$ at $x = \pm L$ and at $y = \pm
L'$. 
Also, from the matching conditions for $v_j^{\pm}$, we have
$$
v^+_x(x,0) - v^-_x(x,0) = 2 \alpha w^+_x(x,0).
$$
Since $\|D^2 w^+_j\|$ is bounded, and
$$ 
\int_{0}^{L'} (w^+_j(\pm L,y))^2 dy \leq C \alpha^2 a^3, \quad
\quad \int_{-L}^{L} (w^+_j(x,L'))^2 dx = 0,
$$ 
it follows from the trace theorem that $w_x^+(x,0)$ is bounded in
$L^2$. 

With this definition of $\tilde{u}$, and suppressing the subscripts $j$ 
and the superscripts $\pm$, we have
\begin{eqnarray}
 {\mathcal I}(u,v,w) & \geq & \int\left[\left(\tilde{u}_x + \frac{1}{2}w_x^2\right)^2 +
\frac{1}{2} (v_x + \tilde{u}_y +  w_x w_y)^2 + (v_y + \frac{1}{2} w_y^2 )^2\right] dx dy
\nonumber
\end{eqnarray}
Given any $\epsilon > 0$, up to extraction of a further subsequence,
 we get
$$ {\mathcal I} \geq \iint ({\tilde{u}}_x)^2 +
({v}_y)^2 + \frac{1}{2}({v}_x + {\tilde{u}}_y)^2 - C
(\|w^*\|^2_{W^{1,4}} + \epsilon)(\|D\tilde{u}\| + \|Dv\|).
$$ The argument from above shows that $v^+_x - v^-_x = 2 \alpha w_x^+$
is bounded in $L^2$.  Combining this with the previous lemma, we see
that $\|D \tilde{u}\| + \|Dv\|$ is bounded in $L^2$.  

$\tilde{u}$ satisfies the boundary conditions $\tilde{u} = 0$ at $x =
\pm L$, $y = \pm L'$. $v$ satisfies the boundary conditions
$$ 
v^{\pm}(\pm L,y) = y \pm \frac{\delta}{2} + \alpha^2 a
\Psi^{\pm}_{1,2}\left(\frac{y}{a}\right) \quad \quad v^{\pm}(x,\pm L') = \pm L' \pm \frac{\delta}{2}.
$$ By the boundedness of $\|Dv\|$, it follows that $v^{\pm}(x,0)$
exist in the sense of traces, and further are in
$L^2(dx)$. Consequently both $\tilde{u}$ and $v^{\pm}$ are bounded in
$L^2(\partial \mathcal{S}^{\pm})$. Combining this with the boundedness
of $\|D \tilde{u}\| + \|Dv\|$, it follows that a further subsequence
$(\tilde{u}_j,v_j)$ converges weakly to $(\tilde{u}^*,v^*)$ in $H^1
\times H^1$, and consequently $(u_j,v_j) \rightharpoonup (u^*,v^*)$
where
$$
u^* = \tilde{u}^* + x
$$ Since $H^1_0$ is weakly closed, it follows that $u^*$ and
${v^\pm}^*$ satisfy the boundary conditions in the sense of traces
(See \cite{evans2} for a detailed argument). We now observe that the
functional $\mathcal{I}$ is convex in $Du, Dv$ and $D^2 w$, and
consequently weakly lower semi-continuous on $H^1 \times H^1 \times
H^2$. This implies that $(u^*,v^*,w^*)$ is a minimizer for
$\mathcal{I}$.
\end{proof}

The existence of a minimizer opens the door to a direct analysis of
the Euler-Lagrange equations for the energy functional
$\mathcal{I}$. We will pursue this approach elsewhere. In this paper,
we restrict ourselves to arguments that do not involve forming the
first variation of $\mathcal{I}$.

\section{Lower Bound} \label{sec:l_bnd}

In this section, we prove a lower bound for the linearized Elastic
energy $\mathcal{I}$ in Eq.~(\ref{eq:unscaled}), by proving a
corresponding result in terms of the scaled energy $I$ in
Eq.~(\ref{eq:scaled}). In the remainder of this section, we will
mostly work with the half sheet $\mathcal{S}^+$, although the same
arguments also hold on $\mathcal{S}^-$. With this understanding, we
will drop the superscripts $\pm$.

As we show in \cite{ridge1}, it follows from the boundary condition $U = 0$ at $X = \pm 1$, 
that
\begin{eqnarray}
I(U,V,W) & \geq & \int \left[(U_X + W_X^2)^2 + W_{YY}^2 \right]dX dY
\nonumber \\ & \geq & \int_0^{\infty} \left[ \frac{1}{2} \left(
\int_{-1}^1 W_X^2 dX \right)^2 + \int_{-1}^1 W_{YY}^2 dX \right] dY
\label{eq:onlyw}
\end{eqnarray}
So, it suffices to prove a lower bound for the functional
$$
E(W) = \int_0^{\infty} \left[ \frac{1}{2} \left(\int_{-1}^1 W_X^2 dX 
\right)^2 + \int_{-1}^1 W_{YY}^2 dX \right] dY 
$$ 

As in \cite{ridge1}, we let $E_b$ and $E_s$ denote the quantities
\begin{equation}
E_b = \int W_{YY}^2 dX dY, \quad \quad E_s = 
\int_0^{\infty}\frac{1}{2} \left(\int_{-1}^1 W_X^2 dX \right)^2 dY.
\label{eq:lw_bnd_energies}
\end{equation}
which we will call the (lower bounds for) the {\em bending} and
{\em stretching} energies respectively.

For every $X$, we define
$$
\rho(X) = \left[\int_0^{\infty} W_{YY}^2(X,Y) dY\right]^{-1}.
$$
$\rho(X)$ is a ``local'' (in $X$) measure of the bending energy,
and $[\rho(X)]^{-1}$ can be thought of as the bending energy density
in $X$ that is obtained by integrating out the $Y$ dependence. 

What we will see below it that $\rho(X)$ is the natural length scale
associated with the ridge as a function of $X$, {\em viz.}, $W_0(X)
\sim \rho(X)$, and the bending energy density in $Y$ decays
rapidly for $Y/\rho(X) \gg 1$ (See Fig.~\ref{fig:ridge}). 

Before we begin the proof of the lower bound, we prove the following
elementary, but very useful result.

\begin{lemma} $f \in H^2$, $\int (f'')^2 = \rho^{-1} < \infty$, 
$f(0) = f_0$, $f'(0) = \beta$. Then we have
$$ \int_0^Y f^2(\eta) d \eta \geq \max_{Z \leq Y,\,\, \theta \in (0,1)}\,(1 - \theta) \,
\left[Z \left(f_0 + \frac{\beta Z}{2}\right)^2 + \frac{Z^3}{12} \left(
\beta^2 - \frac{Z}{\theta \rho}\right) \right] 
$$
\label{lem:ineq}
\end{lemma}

\begin{proof}
By the Sobolev
Embedding theorem, $f$ is a $C^1$ function
and
$$
f(\xi) = f_0 + \beta \xi  + \int_0^\xi f''(\eta) (\xi - \eta) d \eta.
$$ 
Defining $T(\xi)$ by $ T(\xi) = \int_0^\xi f''(\eta) (\xi - \eta) d
\eta$, we have
$$
|T(\xi)|^2 \leq \int_0^\xi [f''(\eta)]^2 d \eta \int_0^\xi
(\xi - \eta)^2 d \eta \nonumber \leq \frac{\xi^3}{3 \rho} 
$$
and integrating this equation in $\xi$ yields
$$
\int_0^Z |T(\xi)|^2 d \xi \leq \frac{Z^4}{12 \rho}.
$$ If $ \tau(Z) = \int_0^Z f^2(\eta) d \eta $, $\tau$ is nondecreasing
in $Z$. Using $f(Z) = f_0 + \beta Z + T(Z)$ and the elementary
inequality
$$
|a + b|^2 \geq (1 - \theta) |a|^2 - \frac{1 - \theta}{\theta} |b|^2,
$$
for all $0 < \theta < 1$, we see that
\begin{eqnarray}
\tau(Z) & \geq & (1 - \theta) \int_0^Z (f_0 +
\beta \eta)^2 d \eta - \frac{1 - \theta}{\theta} \int_0^Z
T^2(\eta) d \eta \nonumber \\ & \geq & (1 - \theta)\left(f_0^2 Z +
\beta f_0 Z^2 + \frac{\beta^2}{3}Z^3\right) - \frac{1 -
\theta}{\theta} \left[ \frac{Z^4}{12 \rho}\right] \nonumber \\
& \geq & (1 - \theta) \left[ Z \left( f_0 + \frac{\beta Z}{2}\right)^2 + \frac{Z^3}{12} \left( \beta^2 - \frac{Z}{\theta \rho}\right)\right] \nonumber   
\end{eqnarray}
The result follows by observing that $\tau$ is nondecreasing, and
optimizing the choice of $\theta$.
\end{proof}

\begin{theorem} ${\mathcal I}(u,v,w)$ is as defined in Eq.~(\ref{eq:unscaled}). 
For all $u \in H^1, v^{\pm} \in H^1$ and $w^{\pm} \in 
H^2 \cap W^{1,4}$ satisfying
the no-stretch boundary conditions
\begin{eqnarray}
u = x, v^{\pm} = \Psi_{1,2}^{\pm}(y/a), w = \sqrt{2}\Phi_{1,2}^{\pm}(y/a), \mbox{ at } x = \pm L, \nonumber 
\end{eqnarray}
and the matching condition $w_y^+(x,0) = w_y^-(x,0) + 2 \alpha$, we have the lower bound
$$ 
{\mathcal I}(u,v^{\pm},w^{\pm}) \geq \min\left(C_1 \alpha^{7/3} \sigma^{5/3}
L^{1/3}, C_2 \alpha^2 \sigma^2 \frac{L}{a} \right).
$$
\label{thm:lbound}
\end{theorem}

\begin{remark} Note that we do not need all the matching conditions in Eq.~(\ref{eq:match}). 
\end{remark}

\begin{remark} The form of the lower bound gives a crossover scale 
$L_a \sim \sigma^{1/3} L^{2/3} \alpha^{-1/3} = L_y$ for $a$. This
scale is the same as the one we obtained in our earlier rescalings.

We will prove the theorem by doing the cases $a \ll L_a$ and $a \gg
L_a$ separately. The result for $a \ll L_a$ is obtained by proving the
scaled version of the statement, {\em viz.},
$$
I \geq C_1 \quad \quad \mbox{for } A < A^*,
$$ for a constant $A^*$ that will be determined below. The case $A >
A^*$ is much easier, and follows immediately as a corollary.
\end{remark}

\begin{remark} In earlier work \cite{ridge1}, we proved the same scaling 
result for $a \ll L_a$, but with extra assumptions on the behavior of
$v$ and $w$ at $y = 0$. We do not know {\em a priori} that these
assumptions are satisfied for a real crumpled sheet. In this theorem,
we have removed these hypothesis, and this result is directly
applicable to crumpled sheets.
\end{remark}

We now begin our proof of the theorem for $a \ll L_a$. As in
\cite{ridge1}, the idea behind the proof is to show that the
stretching energy $E_s$ can be bounded from below by a negative power
of the bending energy $E_b$, so that the total energy $E_s + E_b$
tends to $+\infty$ as $E_b \rightarrow 0$ and $E_b \rightarrow
\infty$. This ensures the existence of a positive lower bound for $E$
(and consequently also for $I$).

We set $\beta^{\pm}(X) = W^{\pm}_Y(X,0)$. The matching condition
therefore is $\beta^{+}(X) = \beta^-(X) + \sqrt{2}$.  Before we prove the
theorem, we collect a few useful results in the following lemmas.

\begin{lemma} 
$$
\int_{-1}^1 \int_0^Y W_X^2 \, dY dX \geq C \max_{Z \leq Y} 
\left[ Z^3 \int_{-1}^{1} \beta(X)^2 dX -2 Z^4 E_b\right] - C'A^3. 
$$
\label{lem:poincare}
\end{lemma}
\begin{proof}
Using Lemma~\ref{lem:ineq} with $f(\xi) = W(X,\xi)$, taking $\theta =
1/2$ and integrating the result in $X$, we see that
$$
\int_{-1}^1 \int_0^Y W^2 \, dY dX \geq \frac{1}{24} 
\max_{Z \leq Y} \left[ Z^3 \int_{-1}^{1} \beta(X)^2 dX - 2 Z^4 E_b\right].
$$ 
The Poincare inequality now yields, 
$$
\int_{-1}^1 \int_0^Y W_X^2 \,
dY dX \geq C \int_{-1}^1 \int_0^Y W^2 \, dY dX - C' \int_0^Y
\left[W^2(-1,\xi) + W^2(1,\xi)\right] d \xi.
$$
The result follows from the observation
$$ 
\int_0^Y \left[W^2(\pm1,\xi)\right] d \xi \leq A^3
\int_0^{\infty} \Phi_{1,2}^2(\eta) d \eta \leq C' A^3
$$ 
\end{proof}

Our proof is based on demonstrating that a small bending energy $E_b$
will lead to a large stretching energy. For $A < \mu$, this idea is
quantified by the following lemma.

\begin{lemma}  Let $B = \int_{-1}^{1} \beta(X)^2 dX$, and let 
$$
\mu = \frac{1}{B}\left(\frac{2 A E_b}{B}\right)^3
$$ 
There is a constant $\mu^* > 0$ such that, if $\mu < \mu^*$, the
stretching energy $E_s$ and the total energy satisfy lower bounds
$$
E_s \geq \frac{C B^7}{E_b^5}, \quad \quad \mbox{and} \quad \quad E \geq E_0 = (5 C B^7)^{1/6}
$$
\label{lem:essential}
\end{lemma}

\begin{proof}
By Jensen's inequality, we have
$$ E_s = \int_0^{\infty}\frac{1}{2} \left(\int_{-1}^1 W_X^2 dX
\right)^2 dY \geq \frac{1}{2}\int_0^{Y}\left(\int_{-1}^1 W_X^2 dX
\right)^2 dY \nonumber \geq \frac{1}{2Y} \left[\int_0^{Y} \int_{-1}^1
W_X^2 dX dY\right]^2
$$ 
Lemma~\ref{lem:poincare} now implies that
$$ \sqrt{E_s} \geq \max_{Y \in \mathbb{R}} 2^{-1/2}\left[C (B Y^{5/2} - 2 E_b
Y^{7/2}) - C' A^3 Y^{-1/2}\right].
$$
Setting $B Y^{5/2} = 2 Y^{7/2} E_b$,
we deduce that a characteristic scale $\tilde{Y}$ for $Y$ is given by
$$
\tilde{Y} = \frac{B}{2 E_b}.
$$
Rescaling $Y$ in terms of $\tilde{Y}$, we obtain
$$
\sqrt{E_s} \geq \frac{C B^{7/2}}{8
E_b^{5/2}}\left[\left(\frac{Y}{\tilde{Y}}\right)^{5/2} \left(1 -
\frac{Y}{\tilde{Y}}\right) - C'
\mu \sqrt{\frac{\tilde{Y}}{Y}}
\right],  
$$ 
where $\mu$ is as defined above, {\em i.e.}
$$
\mu = \frac{1}{B}\left(\frac{2 A E_b}{B}\right)^3.
$$ Observe that $z^{5/2}(1-z)$ has a positive maximum at $z = 5/7 >
0$. The lower bound for the stretching energy $E_s$ follows by
continuity of the function $z^{5/2} (1-z) - C' \mu z^{-1/2}$ with
respect to $\mu$ at $z = 5/7$.

Minimizing $E_s + E_b$, we see that
$$
E =  E_s + E_b \geq E_0 \equiv (5CB^7)^{1/6}
$$ 

\end{proof}

We can now prove the theorem.
\begin{proof}
Set $A^* = (\mu^*)^{1/3}/(2 E_0)$, where $\mu^*$ and $E_0$ are as in
Lemma~\ref{lem:essential}.

The precise statement we will prove is 
$$
E \geq E_0 \min(1,A^*/A)
$$

If $E_b \geq E_0$, there is nothing to prove. Therefore, we can assume
that $E_b < E_0$. It follows that $E_b < E_0$ for each
half-sheet.

We first consider the case $A \leq A^*$. As we argued before, the
sheet has a well defined tangent vector at $y = 0$, and this gives the
matching condition $\beta_+(X) = \beta_-(X) + \sqrt{2}$. This, along with the
convexity of the map $\beta(X) \mapsto B = \int_{-1}^{1} \beta^2(X)
dX$, implies that $B^+ + B^-$ is minimized when $\beta^+(X) = -
\beta^-(X) = \frac{1}{\sqrt{2}}$. Therefore, W.L.O.G. $B^+ \geq
1$. For the half sheet the half-sheet $y \geq 0$, we have --

$A \leq A^*, E < E_0$ and $B \geq 2$ implies that
$$
\mu = 
\frac{1}{B}\left(\frac{2 E_b A}{B}\right)^3 < \mu^*.
$$ 
Lemma~\ref{lem:essential} now implies that $E \geq E_0$, for the
half-sheet $y \geq 0$, and consequently for the whole sheet.

If $A \geq A^*$, we still obtain the conclusion $E \geq E_0$ by the
preceding argument if $E_b$ is so small that
$$
\mu = \frac{1}{B}\left(\frac{2 E_b A}{B}\right)^3 < \mu^*.
$$ 
Therefore, we only need to consider the case $\mu > \mu^*$. W.L.O.G
$B (= B^+) \geq 1$, so that $\mu > \mu^*$ implies
$$
E_b > \frac{(\mu^*)^{1/3}}{2A}
$$
and this gives the desired conclusion.
\end{proof}


\section{Upper bounds} \label{sec:u_bnd}

Our goal is to obtain upper bounds for the functional $I$ that scale
in the same way as the lower bound from the previous section as a
function of the nondimensional parameters in the problem, {\em
viz}. $\epsilon$, $\alpha$ and $A$. This will show that we have
captured the optimal scaling behavior of the elastic energy for a
single ridge in a crumpled sheet.

In particular, we want an upper bound that is a constant (independent
of $\epsilon, \alpha$ and $A$) if $A < A^*$, and an upper bound that
scales as $1/A$ for $A > A^*$. Also, we want upper bounds
that are independent of $\epsilon$ and $\alpha$.

The existence of such upper bounds can be motivated as follows. The
energy in the rescaled variables is given by
\begin{eqnarray}
I(U,V,W) & = & \iint \left[(U_X + W_X^2)^2 +
\frac{\epsilon^{-2/3}}{2} (V_X + U_Y + 2 W_X W_Y)^2 +
\epsilon^{-4/3}(V_Y + W_Y^2)^2 \right] \nonumber
\\ 
& & + \left[ W_{YY}^2 + 2 \epsilon^{2/3} W_{XY}^2 +\epsilon^{4/3}
W_{XX}^2 \right]dX dY. \nonumber
\end{eqnarray}
We would like to show the existence of $(U,V,W)$ satisfying the
boundary conditions such that $I(U,V,W) \leq C(A) < \infty$ uniformly
in $\epsilon$ for $A > 0$. The idea behind the construction of an
appropriate $(U,V,W)$ is as follows. We first pick a smooth $W$
satisfying all the boundary conditions. For this $W$, we will pick an
$V$ such that $V_Y = - W_Y^2$. This equation can (we hope) be solved
for every $X$, along with the appropriate boundary conditions $V(X,Y)
\rightarrow 0$ as $Y \rightarrow \pm \infty$. Once we have $V$, we
determine $U$ by $U_Y = - V_X - 2 W_X W_Y$, again with the appropriate
initial condition for $U$. With such a choice for $U,V$ and $W$, the
energy becomes
$$
\iint \left[(U_X + W_X^2)^2 +  W_{YY}^2 + 2 \epsilon^{2/3} 
W_{XY}^2 +\epsilon^{4/3} W_{XX}^2 \right]dX dY,
$$ and since $U$, $V$ and $W$ are assumed smooth, it easily follows
that there is a finite upper bound, uniform in $\epsilon$ as $\epsilon
\rightarrow 0$. Of course, we are not guaranteed that we have the
right dependence on $A$. Also, we are not guaranteed to get the right
asymptotic shifts in $V^{\pm}(X,.)$.

In the remainder of this section, we will deduce the upper bound by
using ideas similar to the simple argument from above to explicitly
construct smooth functions $(u,v^{\pm},w^{\pm})$ satisfying all the
boundary conditions. With these functions, we can show
$$ {\mathcal I}(u,v,w) \leq \min\left[C \sigma^{5/3} L^{1/3}
\alpha^{7/3} + C' \sigma^2 \alpha^2
\log^+\left(\frac{a}{\sigma}\right), C' \sigma^2 \alpha^2
\frac{L}{a}\right],
$$ where $\log^+ x = \max(\log x, 0)$. This is not exactly the scaling
that we obtained for the lower bounds. In particular, the upper bound
indicates that we are missing some of the relevant physics in our
lower bound if $a \ll \sigma \exp(-\epsilon^{-1/3})$, {\em i.e} $A
\ll \epsilon^{1/3} \exp(\epsilon^{-1/3})$. 

\subsection{Self similar test solutions}

In this section, we show that, for {\em identical, no-stretch}
boundary conditions at $x = \pm L$, with {\em zero asymptotic shift}
($\Delta = 0$), we can construct ``self-similar'' test solutions that
yield the ``correct'' upper bound.

Let $\Phi^{\pm}$ and $\Psi^\pm$ be as in the definition of the
boundary conditions, so that $\Phi^{\pm}(\eta), \Psi^\pm(\eta)$ are
smooth functions, that are supported in $|\eta| \leq K$.

We will choose $W$ and $V$ in the following ``self-similar'' form
$$ 
W^{\pm}(X,Y) = \rho(X) \Phi^{\pm}\left(\frac{Y}{\rho(X)}\right),
\quad \quad V^{\pm}(X,Y) = \rho(X)
\Psi^{\pm}\left(\frac{Y}{\rho(X)}\right),
$$ where $\rho(X)$ is a function smooth that will be chosen later
satisfying $\rho(X) > 0$ for all $X$. The boundary conditions at $X =
\pm 1$ require that $\rho(1) = \rho(-1) = A$.

Let $\eta$ denote $Y/\rho(X)$,
so that $W^{\pm}(X,Y) = \rho(X) \Phi^{\pm}(\eta)$. From this we obtain.
$$
W^{\pm}_{Y}(X,Y) = {\Phi^{\pm}}'(\eta), \quad \quad V^{\pm}_{Y}(X,Y) = {\Psi^{\pm}}'(\eta)
$$ From the matching conditions on $\Phi^{\pm}$ and $\Psi^{\pm}$ at
$\eta = 0$, it is clear that the functions $V^{\pm}$ and $W^{\pm}$
from above satisfy the appropriate matching conditions
(\ref{eq:lin_match}) at $Y = 0$, {\em only if} $\Delta = 0$.

Differentiating in $X$ we obtain,
$$ W^{\pm}_{X} = \rho'(X) \left[ \Phi^{\pm}(\eta) - \eta
{\Phi^{\pm}}'(\eta)\right], \quad \quad V^{\pm}_{X} = \rho'(X) \left[
\Psi^{\pm}(\eta) - \eta {\Psi^{\pm}}'(\eta)\right].
$$
Differentiating once more, we get
\begin{eqnarray}
W^{\pm}_{YY}(X,Y) & = & \frac{1}{\rho(X)} {\Phi^{\pm}}''(\eta), \nonumber \\
W^{\pm}_{XY}(X,Y) & = & - \frac{\rho'(X)}{\rho(X)} \left[\eta {\Phi^{\pm}}''(\eta)\right],
\nonumber \\ 
W^{\pm}_{XX}(X,Y) & = & \frac{1}{\rho(X)} \left( \rho(X) \rho''(X) \left[
\Phi^{\pm}(\eta) -
\eta {\Phi^{\pm}}'(\eta)\right] + \left[\rho'(X)\right]^2 \eta^2 {\Phi^{\pm}}''(\eta)
\right). \nonumber 
\end{eqnarray}
From the no-stretch boundary condition ${\Psi^{\pm}}' + ({\Phi^{\pm}}')^2 = 0$, it
follows that $V^{\pm}_{Y} + [W^{\pm}_{Y}]^2 \equiv 0$.

We would also like $U^{\pm}_{Y} = - V^{\pm}_{X} - 2 W^{\pm}_{X} W^{\pm}_{Y}$. From the
above scalings, we see that
$$ 
V^{\pm}_{X} + 2 W^{\pm}_{X} W^{\pm}_{Y} = \rho'(X) \left[
\Psi^{\pm}(\eta) - \eta {\Psi^{\pm}}'(\eta) + 2
{\Phi^{\pm}}'(\eta)\left(\Phi^{\pm}(\eta) - \eta
{\Phi^{\pm}}'(\eta)\right)\right].
$$
In order that $U^{\pm}_{Y}$ have this scaling behavior, we will set
$$
U^{\pm}(X,Y) = \rho'(X) \rho(X) \Xi^{\pm}\left(\frac{Y}{\rho(X)}\right).
$$
If we choose $\Xi$ such that
$$ {\Xi^{\pm}}' = - \left[ \Psi^{\pm}(\eta) - \eta {\Psi^{\pm}}'(\eta) +
2 {\Phi^{\pm}}'(\eta)\left(\Phi^{\pm}(\eta) - \eta
{\Phi^{\pm}}'(\eta)\right)\right]
$$ Then, we will have $U^{\pm}_{Y} = - V^{\pm}_{X} - 2 W^{\pm}_{X}
W^{\pm}_{Y}$.  

In this approach, we have first order ODEs for $\Xi^{\pm}$, with
boundary conditions $\Xi^{\pm}(\eta) \rightarrow 0$ as $\eta
\rightarrow \pm \infty$, and these ODEs can be solved (in principle)
to yield the functions $\Xi^{\pm}$. These functions are also
required to satisfy the matching condition $U^{+}(X,0) = U^{-}(X,0)$
for all $X$, {\em i.e}, the condition $\Xi^{+}(0) = \Xi^{-}(0)$. 

Since 
the unique solutions for the ODEs determining $\Xi^{\pm}$ are also
supported in $|\eta| \leq K$.

Using the no-stretch condition ${\Psi^{\pm}}' + ({\Phi^{\pm}}')^2 = 0$
in the ODEs for $\Xi^{\pm}$, we get
\begin{eqnarray}
{\Xi^{\pm}}' & = & - \left[ \Psi^{\pm}(\eta) + \eta {\Psi^{\pm}}'(\eta) +
2 {\Phi^{\pm}}'(\eta)\Phi^{\pm}(\eta) \right] \nonumber \\
& = & -\left[\eta \Psi^{\pm}(\eta) + [\Phi^{\pm}(\eta)]^2 \right]'.  \nonumber
\end{eqnarray}
Integrating the above equation, using the fact that $\Phi^{\pm}$ and
$\Psi^{\pm}$ are supported in $|\eta| \leq K$, we deduce that
$\Xi^{\pm}$ is supported in $[-K,K]$ and $\Xi^{\pm}(0) =
[\Phi^{\pm}(0)]^2$.  The matching condition for $\Phi^{\pm}$ now
implies that $\Xi^{+}(0) =\Xi^{-}(0)$. Consequently,
$$
U^{\pm}(X,Y) = \rho'(X) \rho(X) \Xi^{\pm}\left(\frac{Y}{\rho(X)}\right).
$$ satisfies the boundary conditions $U^{\pm} \rightarrow 0$ as $Y
\rightarrow \pm \infty$, and the matching conditions $U^+(X,0) =
U^-(X,0) = 0$. We also require that $U = 0$ at $X = \pm 1$. since
$\rho(\pm1) = A$, we will now require that $\rho'(\pm 1) = 0$.


We will henceforth restrict ourselves to considering the half-sheet $Y
\geq 0$, since the same arguments will also apply to the half-sheet $Y
\leq 0$, and we can now drop the subscripts $\pm$. The above procedure
yields an appropriate test configuration $(U,V,W)$ for boundary
conditions that are identical at $x = \pm L$ and satisfy $\Delta =
0$. We will henceforth refer to this situation as the {\em
self-similar case}. In the remainder of this section, we will consider
this special case, and we will consider the general case in
Sec.~\ref{sec:non_slf_smlr}.

From the above arguments we see that, for appropriate boundary
conditions, it is indeed possible to choose $(U,V,W)$ in the {\em
self-similar form}
\begin{eqnarray}
W(X,Y) & = & \rho(X) \Phi\left(\frac{Y}{\rho(X)}\right), \nonumber \\
V(X,Y) & = & \rho(X) \Psi\left(\frac{Y}{\rho(X)}\right), \nonumber \\
U(X,Y) & = & \rho'(X) \rho(X) \Xi\left(\frac{Y}{\rho(X)}\right). \nonumber
\end{eqnarray}
such that $V_Y + W_Y^2 = 0, U_Y + V_X + 2 W_X W_Y = 0$. With these
choices, we have
$$
I(U,V.W) = \iint \left[(U_X + W_X^2)^2 +  W_{YY}^2 + 2 \epsilon^{2/3} 
W_{XY}^2 +\epsilon^{4/3} W_{XX}^2 \right]dX dY.
$$ A straightforward calculation allows us to estimate the various
terms in this expression. We obtain,
\begin{eqnarray}
\iint (U_X + W_X^2)^2 dX dY & \leq & C \int \left([\rho''(X)
\rho(X)]^2 + [\rho'(X)]^4\right) \rho(X) dX. \nonumber \\ 
\iint W_{YY}^2 dX dY & \leq & C \int\frac{1}{\rho(X)} dX \nonumber \\ 
\iint W_{XY}^2 dX dY & \leq & C \int\frac{[\rho'(X)]^2}{\rho(X)} dX
\nonumber \\ 
\iint W_{XX}^2 dX dY & \leq & C \int \left([\rho''(X)]^2 \rho(X) 
+ \frac{[\rho'(X)]^4}{\rho(X)}\right) dX \label{eq:estimates}
\end{eqnarray}
where $C$ is a constant that only depends on $\Phi$, and $\Psi$.

\subsection{Construction of the upper bound}

In order to prove the claimed upper bound, we need to show the
existence of a smooth $\rho(X)$ such that $\rho(\pm 1) = A$,
$\rho'(\pm 1) = 0$, and all the terms in $I$ are bounded uniformly in
$\epsilon$ as $\epsilon \rightarrow 0$ with the appropriate dependence
$A$.

Lobkovsky's \cite{lobkovsky} analysis motivates the choice $\rho(X)
\sim (1 - X)^{2/3}$ near $X = 1$. For this choice however, the
contribution of $W_{XX}^2$ is given by
$$
[\rho''(X)]^2 \rho(X) \sim \frac{[\rho'(X)]^4}{\rho(X)} \sim (1 - X)^{-2},
$$ and is not integrable near $X = 1$. Therefore the choice $\rho(X)
\sim (1 - |X|)^{2/3}$ will not yield an upper bound for $I$.

In our analysis of the lower bound, we ignored the contribution of
$W_{XX}$ to the energy, and obtained results that agree with
Lobkovsky's boundary layer analysis. This suggests that the scaling
$\rho(X) \sim (1 - |X|)^{2/3}$ might still be appropriate in regions
where the contribution of the $W_{XX}$ and the $W_{XY}$ terms are
small. However, close to the boundaries near $X = \pm 1$, the dominant
energy balance is different, and we need to modify the behavior of
$\rho$ to account for this.

For small $X$, we expect that $\epsilon^{4/3} W_{XX}^2 \sim W_{YY}^2$
is the leading order balance for the energy in $I$. Since $Y \sim
\rho(X)$, it follows that $\rho(X) \sim \epsilon^{-1/3}(1- X)$ near $X
= 1$, and similarly $\rho(X) \sim \epsilon^{-1/2}(1+X)$ near $X =
-1$. Note that these behaviors match $\rho(X) \sim (1-X)^{2/3}$
(respectively $(1+X)^{2/3})$) when $1 - X \sim \epsilon$ (respectively
$1+X \sim \epsilon$).

We also have the boundary conditions $\rho(\pm 1) = A$ and $\rho'(\pm
1) = 0$. This suggests $\rho(X) \approx A$ for $(1 - |X|) < A
\epsilon^{1/3}$ if $A \ll \epsilon^{2/3}$. Therefore we will choose
$\rho(X)$ with the following behavior.
\begin{itemize}
\item In the case $A \ll \epsilon^{2/3} \ll 1$, 
$$ \rho(X) \sim \left\{ \begin{array} {cc} A & (1 - |X|) \lesssim
A\epsilon^{1/3} \\ 
\epsilon^{-1/3} (1 - |X|) & A\epsilon^{1/3} \ll (1 -|X|) \ll
\epsilon \\ (1 - X^2)^{2/3} & \epsilon \ll (1-|X|) \sim 1 \end{array}
\right.
$$
\item In the case $\epsilon^{2/3} \lesssim A \ll 1$ 
$$
\rho(X) \sim \left\{ \begin{array} {cc} A & (1 - |X|) \lesssim A \\
(1 - X^2)^{2/3} & \epsilon \ll (1-|X|) \sim 1 \end{array} \right.
$$
\item In the case $A \sim 1$, we set $\rho(X) = A$.
\end{itemize}

We will now make the above considerations precise. Let $\varphi \in
C_c^{\infty}$ be a smooth, nonnegative, non-increasing function, that
is identically one on $(- \infty,1/2]$ and zero on $[2,\infty)$.  An
example of such a function is illustrated in Fig. Also,
$\bar{\varphi}$ will denote the complementary function $1 - \varphi$.

\begin{lemma} For $A > 0$ and $\epsilon > 0$, let 
\begin{eqnarray}
g(z) & = & A \varphi \left(\frac{z}{A \epsilon^{1/3}}\right) +
z \epsilon^{-1/3}  \bar{\varphi}\left(\frac{z}{A
\epsilon^{1/3}}\right) \varphi\left(\frac{z}{\epsilon}\right) + z^{2/3}
\bar{\varphi}\left(\frac{z}{\epsilon}\right), \nonumber \\
h(z) & = &
z \epsilon^{-1/3} \chi_{[\frac{1}{2} A \epsilon^{1/3}, 2\epsilon]} + z^{2/3}
\chi_{[\frac{1}{2}\epsilon,\infty)}, \nonumber
\end{eqnarray} 
where $\chi$ denotes the characteristic function. Then $g$ is a smooth
function. Also, there exist constants $c$ and $C$ such that $g$
satisfies the following inequalities $\forall z$
\begin{eqnarray}
g(z) & \geq & c \left[A \chi_{[0,2 A \epsilon^{1/3}]} + h(z)\right]
\nonumber \\ 
g(z) & \leq & C \left[A \chi_{[0,2 A \epsilon^{1/3}]} +
h(z) \right] \nonumber \\ 
g'(z) & \leq & C h(z)
\left[1 + \frac{1}{z}\right]  \nonumber \\ 
g''(z) & \leq & C h(z) \left[1 + \frac{1}{z^2}\right] \nonumber
\end{eqnarray}
\label{lem:crossover}
\end{lemma}
\begin{proof}
We begin with an elementary observation. Near $z \sim A
\epsilon^{1/3}$, the functions $g = A$ and $g = z \epsilon^{-1/3}$ are
comparable, {\em i.e}, there exist constants $\theta,\Theta$ {\em independent of
$\epsilon,A$}, such that $\theta A \leq z \epsilon^{-1/3} \leq \Theta A$ for
$\frac{1}{2} A \epsilon^{1/3} \leq z \leq 2 A
\epsilon^{1/3}$. Similarly, near $z \sim \epsilon$, we have $\theta z
\epsilon^{-1/3} \leq z^{2/3} \leq \Theta z \epsilon^{-1/3}$ for
$\frac{1}{2} \epsilon \leq 2 \epsilon$.

This observation implies that $g(z) \geq c \left[A \chi_{[0,2 A
\epsilon^{1/3}]} + h(z)\right]$. The inequality $g(z) \leq C \left[A
\chi_{[0,2 A \epsilon^{1/3}]} + h(z)\right]$ is elementary and follows
from the boundedness of $\varphi$.

We also observe that, for all $l > 0$ and $n = 1,2,3,\ldots$,
$$ 
\left|\frac{d^n}{dz^n} \varphi\left(\frac{z}{l}\right) \right|\leq \frac{C_n}{l^n}
\chi_{[\frac{1}{2}l, 2l]} \leq \frac{C_n}{z^n} \chi_{[\frac{1}{2}l,
2l]},
$$ and the same inequality also holds for the complementary function
$\bar{\varphi}$.  We also have the elementary inequality 
$$ \left|\frac{d^n}{dz^n} z^{\alpha} \right| \leq \frac{C_{n,\alpha}
z^{\alpha}}{z^n}, \quad \quad z > 0
$$ Differentiating $g$ twice, using all of the above observations, and
recognizing that $z^{-1} \leq 1 + z^{-2}$, we get the inequalities
\begin{eqnarray}
|g'(z)| & \leq & C h(z)
\left[1 + \frac{1}{z}\right]  \nonumber \\ 
|g''(z)| & \leq & C h(z) \left[1 + \frac{1}{z^2}\right] \nonumber
\end{eqnarray}
This proves the lemma.
\end{proof}

\begin{lemma} 
For boundary conditions that support a self-similar 
test function, we have an upper bound
$$ I(U^*,V^*,W^*) \leq C \left[ 1 +
\epsilon^{1/3}\log^{+}\left(\frac{\epsilon^{2/3}}{A}\right) \right].
$$
\label{lem:slf_smlr}
\end{lemma}

\begin{proof}
We set $\rho(X) = g(1-X^2)$, where $g$ is as defined in
lemma~\ref{lem:crossover}. Note that $g(0) = A, g'(0) = 0$ implies
$\rho(\pm1) = A, \rho'(\pm 1) = 0$. 

We begin with a few observations.
\begin{enumerate}
\item Our definition makes 
$\rho$ an even function of $X$, and by lemma~\ref{lem:crossover},
$\rho$ is smooth. 
\item $z'(X) = - 2 X$ and $z''(X) =
-2$ are bounded for $X \in [-1,1]$, so that the inequalities in
lemma~\ref{lem:crossover} also hold in terms of $\rho$ with the
understanding $z = 1 - X^2$.
\item If $A < \epsilon^{2/3} < 1$, the sets $[0,A\epsilon^{1/3}), [A
\epsilon^{1/3}, \epsilon)$ and $[\epsilon,1]$ give a disjoint
partition of $[0,1]$. If $A > \epsilon^{2/3}$, the sets
$[0,A\epsilon^{1/3})$, and  $[\epsilon,1]$ cover $[0,1]$.
\end{enumerate}

We can now estimate the various terms in (\ref{eq:estimates}) using
the results of lemma~\ref{lem:crossover}, and the observations from
above.

We first estimate the $(U_X+W_X^2)^2$ term. Using $(1+z^{-1})^4 \leq C
(1+z^{-2})^2 \leq C(1+z^{-4})$, we obtain
\begin{eqnarray}
\int \left([\rho''(X) \rho(X)]^2 + [\rho'(X)]^4\right) \rho(X) dX &
\leq & C \int_0^1 [h(z)]^5 \left(1 + \frac{1}{z^4} \right) dz,
\nonumber \\ & \leq & C \epsilon^{-5/3}
\int_{\frac{1}{2}A\epsilon^{1/3}}^{2 \epsilon} z^5 \left(1 +
\frac{1}{z^4} \right) dz \nonumber \\ & & + C \int_{\frac{1}{2}
\epsilon}^{1} z^{10/3} \left(1 + \frac{1}{z^4} \right) dz,\nonumber \\
& \leq & C\left[ \epsilon^{1/3} + 1\right] \nonumber
\end{eqnarray}
The $W_{YY}^2$ term yields
\begin{eqnarray}
\int\frac{1}{\rho(X)} dX & \leq & C \int_0^{2 A\epsilon^{1/3}}
\frac{dz}{A} + C \int_{\frac{1}{2} A\epsilon^{1/3}}^{2 \epsilon}
\frac{dz}{z \epsilon^{-1/3}} + C \int_{\frac{1}{2}\epsilon}^{1}
\frac{dz}{z^{2/3}} \nonumber \\ & \leq & C \left[\epsilon^{1/3}
\left(\log\left(\frac{\epsilon^{2/3}}{A}\right)+1\right) + 1\right] \nonumber
\end{eqnarray}
We will now estimate the $W_{XX}$ term, since a bound for the $W_{XY}$
term can be obtained from the bounds for the $W_{YY}$ and the $W_{XX}$
terms. We have
\begin{eqnarray}
\int \left([\rho''(X)]^2 \rho(X) + \frac{[\rho'(X)]^4}{\rho(X)}\right)
dX & \leq & C \int_0^1 [h(z)]^3 \left(1 + \frac{1}{z^4}\right) dz \nonumber \\
& \leq & C \epsilon^{-1}
\int_{\frac{1}{2}A\epsilon^{1/3}}^{2 \epsilon} z^3 \left(1 +
\frac{1}{z^4} \right) dz \nonumber \\ & & + C \int_{\frac{1}{2}
\epsilon}^{1} z^2 \left(1 + \frac{1}{z^4} \right) dz,\nonumber \\
& \leq & \frac{C}{\epsilon}  \left(\log\left(\frac{\epsilon^{2/3}}{A}\right)+1\right) \nonumber
\end{eqnarray}
Using 
$$ \frac{[\rho'(X)]^4}{\rho^2(X)} = \frac{[\rho'(X)]^4}{\rho(X)} \cdot \frac{1}{\rho(X)},
$$
and the Cauchy-Schwarz inequality, the above estimates yield
$$
\int  \frac{[\rho'(X)]^2}{\rho(X)} 
dX \leq C \left[\epsilon^{-1/3}
\left(\log\left(\frac{\epsilon^{2/3}}{A}\right)+1\right) + \epsilon^{-2/3} \right]
$$ We have thus bounded the $W_{XY}$ term. Using these estimates in
$$
I(U,V.W) = \iint \left[(U_X + W_X^2)^2 +  W_{YY}^2 + 2 \epsilon^{2/3} 
W_{XY}^2 +\epsilon^{4/3} W_{XX}^2 \right]dX dY,
$$ we see that, for $A < \epsilon^{2/3}$ the self-similar test
function has
$$
I \leq C \left[ 1 + \log\left(\frac{\epsilon^{2/3}}{A}\right) \right].
$$ If $\epsilon^{2/3}<A$, the sets $[0,A \epsilon^{1/3})$ and
$[\epsilon,1]$ cover $[0,1]$. In this case, the same analysis as above
gives $I \leq C$. Combining this with the above result proves the
lemma.
\end{proof}
Setting $\rho(X) = A$, a direct calculation using
(\ref{eq:estimates}), shows that $I \leq C/A$. Combining this result
with the lemma~\ref{lem:slf_smlr}, and ``unscaling'' these results, we get
\begin{theorem}
If the no stretch boundary conditions are identical at $x = \pm L$,
and have $\Delta = 0$, we have the upper bound
$$ I(U^*,V^*,W^*) \leq C\min\left( \left[ 1 + \epsilon^{1/3}
\log^{+}\left(\frac{\epsilon^{2/3}}{A}\right)
\right],\frac{C'}{A}\right).
$$
\label{thm:slf_smlr}
\end{theorem}

\subsection{Non-self similar test solutions} \label{sec:non_slf_smlr}

We will now consider upper bounds for the general case, {\em i.e}, for
the situation where the boundary conditions at $x = \pm L$ are not
identical and/or $\Delta \neq 0$.

\begin{figure}[htbp]
  \begin{center}
\centerline{\epsfig{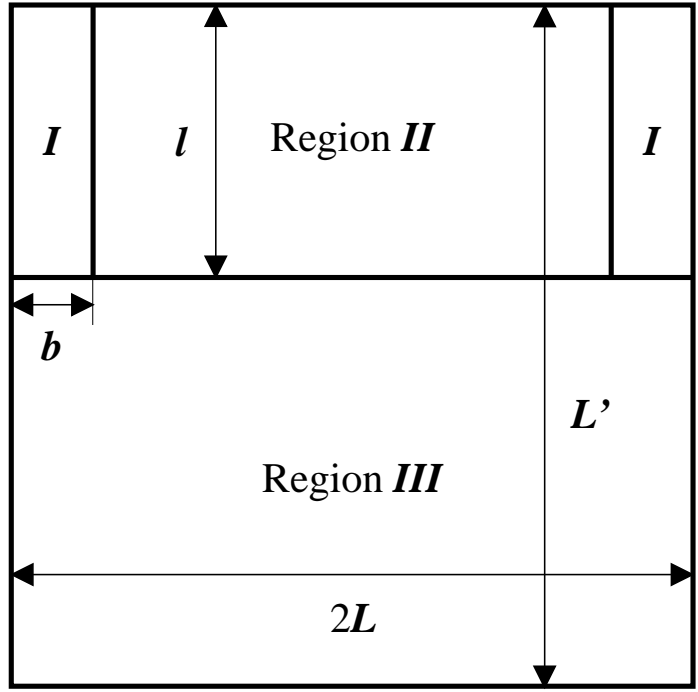}}
\end{center}
\caption{A schematic illustration of the various regions, and their
sizes.}
\label{fig:matched}
\end{figure}

The strategy of the proof will be the following:
\begin{enumerate}
\item We introduce boundary layers (in $x$) near the boundaries $x =
\pm L$ of width $x = b$. In these boundary layers, we connect the
boundary condition at $x = \pm L$ with identical profiles that have
$\Delta = 0$, at a distance $a$ from the boundaries.
\item In the region, $L-|x|\geq b, |y| \leq K L_y$, we use our {\em
self similar construction} from the last section.
\item In the remaining region, we introduce a {\em small, uniform}
strain $\gamma_{yy} \sim \delta/L'$, to get the appropriate asymptotic
shift $\delta$.
\end{enumerate}
The various regions are illustrated in Fig.~\ref{fig:matched}. The
idea behind this construction is from {\em matched asymptotic
expansions}. The self-similar solutions from the last section play the
role of the {\em outer solutions} in $x$, but the {\em inner
solutions} in $y$!

Before we rigorously construct a test solution that gives the upper
bound, we present a heuristic scaling argument that motivates our
choices for the length scales $b$ and $l$ in our test solution. 

The boundary conditions at $x = \pm L$ are given by no-stretch
profiles, whose (curvature) length scale is $a$. We assume that we can
find a no-stretch profile with zero asymptotic shift, and curvature
$a^{-1}$. This will be the profile of the test solution at $x = \pm (L
- b)$. Also, we set $u = x$ in these layers, so that $u_x = 1, u_y =
0$.

In the boundary layer(s) $L-b \leq |x| \leq L$, $w(x,y)$ is supported
in $|y| \leq K a$. Also, $w \sim O(\alpha a)$ so that
$$
w_x \sim \frac{\alpha a}{b}, \quad \quad w_{xx} \sim \frac{\alpha a}{b^2}, 
w_{yy} \sim \frac{\alpha}{a},
$$ and all of these derivatives are supported in $|y| \leq K a$. This
gives $\gamma_{xx} \sim w_x^2 \sim \alpha^2 a^2/b^2$

If we assume that there is no stretching in the $y$ direction, so that
$\gamma_{yy} = 0$, we see that $v_y \sim 1 - \frac{1}{2}w_y^2$. Consequently, we
see that the asymptotic shift $\delta$ is a function of $x$, and by
integrating $w_y^2$ we see that $\delta(x) \sim O(\alpha^2 a)$, and
this gives $v_x \sim \delta_x \sim \alpha^2 a/b$. Also, since this
$v_x$ arises from the difference in the asymptotic shifts, it is
supported in $|y| \leq l$, and not only in $|y| \leq a$, as for
$w_y$. This gives $\gamma_{xy} \sim v_x \sim \alpha^2 a/b$.

Putting all of this together, the elastic energy in the two boundary
layers  comprising Region I is
\begin{eqnarray}
\mathcal{E}_I & \sim & \gamma_{xx}^2 a b + \gamma_{xy}^2 l b + 
\sigma^2\left(w_{xx}^2 + w_{yy}^2 \right) a b \nonumber \\
& \sim & \frac{\alpha^4 a^5}{b^3} + \frac{\alpha^4 a^2 l}{b} + 
\frac{\sigma^2 \alpha^2 b}{a} + \frac{\sigma^2 \alpha^2 a^3}{b^3} \nonumber
\end{eqnarray}

In Region II, we use the self similar construction from above, {\em
 i.e.}, the dominant energy balance is between the curvature $w_{yy}$
 and the strain $\gamma_{xx}$. Following the scaling argument in
 Ref.~\cite{science.paper}, we get
$$
\mathcal{E}_{II} \sim \frac{\alpha^4 l^5}{L^3} + \frac{\alpha^2 \sigma^2 L}{l}
$$

In region III, we are connecting $v = y$ (zero asymptotic shift) with
$v = y \pm \delta$. We can do this with a profile that has a uniform
strain $\gamma_{yy} \sim \delta/(L'-l)$. Consequently,
$$ 
\mathcal{E}_{III} \sim \left(\frac{\delta}{L'-l}\right)^2 L (L'-l) \sim
\frac{\alpha^4 a^2 L}{L'-l}.
$$ Note that $\mathcal{E}_{III}$ is completely determined by the
boundary conditions, and is {\em independent} of any choice we make
for $b$ and $l$ provided that $l \ll L'$.

For any given $a,\sigma,L,L'$, $\mathcal{E} = \mathcal{E}_{I} +
  \mathcal{E}_{II} + \mathcal{E}_{III}$ diverges as $b,l \rightarrow 0,
  \infty$ independently. Consequently there are optimal finite,
  nonzero choice $b = b^*$ and $l = l^*$.

There are two {\em scaling regimes} of interest. In physically
realistic situations, $\sigma \lesssim a$. If $a = C \sigma$, with $C$
staying $O(1)$ as $\sigma \rightarrow 0$, optimizing $l$ and $b$ gives
$$
l^* \sim \alpha^{-1/3} \sigma^{1/3} L^{2/3}, \quad \quad b^* \sim \alpha^{5/6}
 \epsilon^{2/3} L^{1/3}.
$$
The energies in the three regions are 
$$
\mathcal{E}_{I} \sim \alpha^{17/6} \sigma^{5/3} L^{1/3}, \quad \quad 
\mathcal{E}_{II} \sim \alpha^{7/3} \sigma^{5/3} L^{1/3}, \quad \quad 
\mathcal{E}_{III} \sim \alpha^4 \sigma^2  \frac{L}{L'}. 
$$ For $\sigma \ll L,L'$, the energy $\mathcal{E}_{III}$ is asymptotically
negligible. $\mathcal{E}_{I} \sim \alpha^{1/2} \mathcal{E}_{II}$, so the
two energies scale in the same way for $\alpha \sim O(1)$, but the
energy in the boundary layers is asymptotically negligible in the
small angle limit.

We can also consider the situation $a \ll \sigma$. In this case, the
natural scaling regime is $a \sim C \sigma^{5/3} L^{-2/3}$, with $C$
staying $O(1)$ as $\sigma \rightarrow 0$. Optimizing $l$ and $b$ gives
$$ l^* \sim \alpha^{-1/3} \sigma^{1/3} L^{2/3}, \quad \quad b^* \sim
\epsilon^{5/3} L^{-2/3} \sim a.
$$
The energies in the various regions are
$$
\mathcal{E}_{I} \sim \alpha^{2} \sigma^{2}, \quad \quad 
\mathcal{E}_{II} \sim \alpha^{7/3} \sigma^{5/3} L^{1/3}, \quad \quad 
\mathcal{E}_{III} \sim \alpha^4 \sigma^2  \left(\frac{\sigma}{L}\right)^{1/3}
\frac{\sigma}{L'}. 
$$ For $\sigma \ll L,L'$, the energies $\mathcal{E}_{I}$ and
$\mathcal{E}_{III}$ are asymptotically negligible, and the energy is
determined, essentially by the the self-similar solution in Region II,
in the limit $\epsilon = (\sigma/L\alpha) \rightarrow 0$.

We will now use these scaling results as motivation, and rigorously
construct test solutions that give the appropriate upper bound in
situations where the boundary conditions are non-identical, or have
nonzero asymptotic shift. Before we begin our construction, we first
show that there exist no-stretch profiles with zero asymptotic shift.

\begin{lemma} $\exists K_0$, such that $\forall K \geq K_0,$ 
there exist smooth $\phi,\psi$ supported on $[0,K]$ such
that
$$
\phi'(0) = 1, \quad \quad \int_0^K [\phi''(\eta)]^2 d \eta = \frac{1}{2},
$$
and
$$
\int_0^{K} \left[\phi'(\eta)\right]^2 d \eta + \phi(0) = 0
$$ 
\label{lem:chopped_off}
\end{lemma}

\begin{proof}
Let $\zeta$ be a smooth nonnegative function such that $0 \leq \zeta
\leq 1$, $\zeta(x) \equiv 1$ for $x \leq -1$ and $\zeta(x) \equiv 0$
for $x \geq 0$. Let $\varpi$ be a smooth function supported in
$[-1,0]$, that is odd about $-\frac{1}{2}$, and is not identically
zero.

We set $\varsigma(x) = \zeta(x) + q \varpi(x)$. The map
$$
q \mapsto \Delta(q) \equiv  \int_{-1}^0 \varsigma(\xi)\left[1 - \varsigma(\xi)\right] d \xi
$$ is clearly continuous. The reason for calling this map $\Delta$
will become clear below. 

We have
$$
\Delta(0) = \int_{-1}^0 \zeta(\xi)\left[1 - \zeta(\xi)\right] d \xi \geq 0
$$ since $0 \leq \zeta \leq 1$. Also,
\begin{eqnarray}
\Delta(q) & = &  \int_{-1}^0 \zeta(\xi) d \xi + q \int_{-1}^0 \varpi(\xi) d \xi 
- \int_{-1}^0 \left[\zeta(\xi) + q \varpi(\xi)\right]^2 d \xi \nonumber \\
& \leq &  \int_{-1}^0 \zeta(\xi) d \xi + 2  \int_{-1}^0 \zeta^2(\xi) d \xi - \frac{q^2}{2}  \int_{-1}^0 \varpi^2(\xi) d \xi \nonumber \\
& \leq & 3 - C q^2 \nonumber
\end{eqnarray}
where we have used $|\zeta| \leq 1$ and $\varpi$ is not identically
zero in passing to the last line. 

Since $\Delta(0) \geq 0$, and $\Delta(q) \rightarrow - \infty$ as $q
\rightarrow - \infty, \exists q^* < \infty$ such
that $\Delta(q^*) = 0$. 

We set $\varsigma = \zeta + q^* \varpi$, and define $l$ by
$$
\int_{-1}^{0} \left[\varsigma'(x)\right]^2 dx = \frac{l}{2}.
$$
Let $K_0 = 2l$. For a given $K$, we will set 
$$
\phi(x) = \int_K^x \varsigma\left(\frac{\xi - K}{l}\right) d\xi.
$$ Then, $\phi'(x) = \varsigma\left(\frac{x - K}{l}\right)$, so that
$\phi'(x) \equiv 1$ for $x \leq l \leq K -l$. Also,
$$
 \int_0^K [\phi''(x)]^2 d x = \frac{1}{l^2}  \int_0^K \left[\varsigma'
\left(\frac{x - K}{l}\right)\right]^2 d x = \frac{1}{2}.
$$
Since $\phi$ is supported on $[0,K]$, it follows that 
$$
\phi(0) = - \int_0^K \phi'(x) dx.
$$
Since $\phi(x) = 1$ for $0 \leq x \leq K -l$, it follows that 
\begin{eqnarray}
\int_0^{K} \left[\phi'(\eta)\right]^2 d \eta + \phi(0) & = &
\int_{K-l}^{K} \left[(\phi'(\eta))^2 - \phi'(\eta)\right] d \eta
\nonumber \\ & & = \Delta(q^*) l = 0. \nonumber 
\end{eqnarray}
We now set 
$$
\psi(x) = \int^K_x[\phi'(\xi)]^2 d \xi.
$$
\end{proof}
We can extend the functions $\phi,\psi$ to a no-stretch profile on
$\mathbb{R}$ by setting
$$ \tilde{\Phi}^+(\eta) = \tilde{\Phi}^-(-\eta) = \phi(\eta), \quad \quad \tilde{\Psi}^+(\eta)
= -\tilde{\Psi}^-(-\eta) = \psi(\eta), \quad \quad \eta \in [0,K]
$$ 
and $\tilde{\Phi}^{\pm} = \tilde{\Psi}^{\pm} = 0$ otherwise. This profile satisfies
the matching conditions and has zero asymptotic shift.

We now have to prove that, by introducing a thin boundary layer near
the boundaries, we can ``connect'' the prescribed boundary conditions
to the profile we constructed in lemma~\ref{lem:chopped_off}, without
incurring a large energy penalty. We begin with the following lemma
which estimates norms of the first derivatives of a no stretch profile
$(\Phi^{\pm},\Psi^{\pm})$ in terms of the curvature of the profile,
and it's support.

\begin{lemma} $(\Phi^{\pm},\Psi^{\pm})$ is a normalized no-stretch 
profile with support $K$, {\em i.e.}, ${\Psi^{\pm}}' + [{\Phi^{\pm}}']^2 =
0$, $\Phi^{\pm}(\eta) = \Psi^{\pm}(\eta) = 0$ if $|\eta| \geq K$, 
$$
{\Phi^-}'(0) = {\Phi^+}'(0) - \sqrt{2}
$$
and 
$$
\int_0^{\infty} \left[{\Phi^+}''(\eta)\right]^2 d \eta +
\int_{-\infty}^0 \left[{\Phi^-}''(\eta)\right]^2 d \eta = 1.
$$
It then follows $K \geq 1/2$ and
\begin{eqnarray}
\frac{1}{K^3}\sup_{0 \leq \eta \leq K} [\Phi^{+}(\eta)]^2 + \,\frac{1}{K^7}
\int_0^{K} [{\Phi^+}(\eta)]^4 d \eta \,+\,
\frac{1}{K^4}\int_0^{K} [{\Phi}^+(\eta)]^2 d \eta & \leq  & C  \nonumber \\
\frac{1}{K^2} \left[\int_0^{K} [{\Phi^+}'(\eta)]^2 d \eta + \sqrt{2} \Phi^+(0)\right]
& \leq & C
\end{eqnarray}
where $C$ is a universal constant (independent of $K$). A similar
result also hold for $\Phi^-$.
\label{lem:bounds}
\end{lemma}
\begin{proof}
This follows immediately from the Sobolev embedding theorem in $\mathbb{R}^1$
\cite{evans,ziemer}. Since ${\Phi^+}'(K) = 0$, we have
\begin{eqnarray}
|{\Phi^+}'(\eta)| & = & \left|\int^K_{\eta} {\Phi^+}''(\xi) d \xi\right|
\nonumber\\ 
& \leq & \left[ \int_\eta^K [{\Phi^+}''(\xi)]^2 d \xi \cdot
  \int_\eta^K 1 d \xi\right]^{1/2} \nonumber \\ 
& \leq & \sqrt{K - \eta} \nonumber
\end{eqnarray}
Using ${\Phi^-}'(0) = {\Phi^+}'(0) - \sqrt{2}$, it follows that $2 \sqrt{K}
\geq \sqrt{2}$, which implies that $K \geq 1/2$.  

The remaining inequalities follows from integrating
$|{\Phi^+}'(\eta)| \leq \sqrt{K - \eta}$ and using $\Phi^+(K) = 0$,
and $K \geq 1/2$ so that $K^2 > \frac{1}{2}K^{3/2}$.
\end{proof}

\begin{lemma} $\Phi^{\pm}_{1,2}$ and $\Psi^{\pm}_{1,2}$ are smooth on 
$\mathbb{R}^{\pm}$ and supported on $[0,\pm K]$. Further, for $i =
1,2$, ${\Psi^{\pm}_i}' + [{\Phi^{\pm}_i}']^2 =0$, and
$$ 
\int_0^{\infty} \left[{\Phi^+}''(\eta)\right]^2 d \eta +
\int_{-\infty}^0 \left[{\Phi^-}''(\eta)\right]^2 d \eta = 1.
$$ Also, $\Phi_{1,2}^{\pm}$ and $\Psi_{1,2}^{\pm}$, satisfy the
matching conditions (\ref{eq:lin_match}), where the asymptotic shifts
of the two profiles are given by
$$ 
-\left[\int_{-\infty}^0 \left[{\Phi^-_i}'(\eta)\right]^2 d \eta + \sqrt{2}\Phi^-_i(0) + 
\int_0^{\infty} \left[{\Phi^+_i}'(\eta)\right]^2 d \eta + \sqrt{2}\Phi^+_i(0)\right] = \Delta_i.
$$ 

$\zeta$ is a nonnegative, nondecreasing $C^{\infty}$ function such that
$\zeta(x) \equiv 0$ for $x \leq 0$ and $\zeta(x) = 1$ for $x \geq
1$. Let
\begin{eqnarray}
U(X,Y) & = & 0\nonumber\\ 
W^{\pm}(X,Y) & = &
\zeta\left(\frac{X}{B}\right)\Phi^{\pm}_2(Y) + \left[1 -
\zeta\left(\frac{X}{B}\right)\right]\Phi^{\pm}_1(Y) \nonumber \\
V^{\pm}(X,Y) & = & - \int_{\pm K}^{Y} \left[ W^{\pm}_Y(X,Y) \right]^2 d Y + Y 
\pm \frac{1}{2}\Delta(X)\nonumber
\end{eqnarray}
where 
$$ 
\Delta(X)= -\left[\int_{-\infty}^0 \left[W^-_Y(X,\eta) \right]^2 d
\eta + \int_0^{\infty} \left[W^+_Y(X,\eta)\right]^2 d \eta + 2
\sqrt{2} W^+(X,0)\right]
$$ Then, $(U,V^{\pm},W^{\pm})$ satisfies the matching conditions at $Y
= 0$, and we have
\begin{eqnarray}
\int_{Y=0}^{M}\int_{X=0}^B (U_X + W_X^2)^2 \, dX dY & \leq &
\frac{C}{B^3} \nonumber \\ \int_{Y=0}^{M}\int_{X=0}^B (V_X + U_Y + 2
W_X W_Y)^2 \, dX dY & \leq & \frac{CM}{B} \nonumber \\
\int_{Y=0}^{M}\int_{X=0}^B (V_Y + W_Y^2)^2 \, dX dY & = & 0 \nonumber \\
\int_{Y=0}^{M}\int_{X=0}^B W_{YY}^2 \,dX dY & \leq & C B \nonumber \\
\int_{Y=0}^{M}\int_{X=0}^B W_{XX}^2 \,dX dY & \leq & \frac{C}{B^3} \nonumber
\end{eqnarray}
where $C$ is a constant that only depends on $K$.
\label{lem:bndry_lyr}
\end{lemma}
\begin{proof}
Most of the inequalities follow by direct integration using the
definition of $(U,V^{\pm},W^{\pm})$, and using
Lemma~\ref{lem:bounds}. The only result that needs proof is
$$
\int_{Y=0}^{M}\int_{X=0}^B (V_X + U_Y + 2
W_X W_Y)^2 \leq \frac{CM}{B}.
$$ Since $W(X,.)$ is supported in $[0,K]$, and $V^{\pm}(X,Y) = Y \pm
\frac{1}{2} \Delta(X)$ for $|Y| \geq K$, we have $V_X = \Delta_X$. From
the definition of $\Delta(X)$, we see that
\begin{eqnarray}
\frac{d}{dX}\Delta(X) & = & - \frac{2}{B}
\zeta'\left(\frac{X}{B}\right) \left[ \int_{-K}^0 W_Y^-(X,\eta)
[\Phi_2^-(\eta) - \Phi_1^-(\eta)] d \eta \right. \nonumber \\ 
& & \left. + \int_{0}^K W_Y^+(X,\eta)
[\Phi_2^+(\eta) - \Phi_1^+(\eta)] d \eta + \sqrt{2}(\Phi_2(0) -
\Phi_1(0))\right] \nonumber
\end{eqnarray}
Estimating the integrals by the Cauchy-Schwarz inequality, and using
Lemma~\ref{lem:bounds}, we see that
$$
|\Delta_X| \leq \frac{C}{B} \left|\zeta'\left(\frac{X}{B}\right)\right| \quad \implies \quad 
\int_0^B [\Delta_X(\xi)]^2 d\xi \leq \frac{C}{B}.
$$
We also have,
\begin{eqnarray}
\frac{d}{dY} (V_X + U_Y + 2 W_X W_Y) & = & V_{XY} + 2 W_{XY} W_Y + 2 W_X W_{YY} \nonumber \\
& = & \frac{d}{dX}(V_Y + W_Y^2)  + 2 W_X W_{YY}  = 2 W_X W_{YY}  \nonumber 
\end{eqnarray}
Using $U \equiv 0$, and $V_X = \Delta_X, W_X = W_Y = 0$ for $Y \geq K$, we get
\begin{eqnarray}
\int_0^{M} (V_X + U_Y + 2 W_X W_Y)^2 d \eta & = & (M - K) \Delta_X^2
+ \int_0^{K} (V_X + U_Y + 2 W_X W_Y)^2 d \eta \nonumber \\
& \leq & C M \Delta_X^2 + C K^2 \int_0^{K}  W_X^2 W_{YY}^2 d \eta. \nonumber
\end{eqnarray}
using the Poincare inequality. By the convexity of the map,
$$
W(X,.) \mapsto \int_0^{K} W_{YY}^2 d \eta,
$$ 
and the normalization of $\Phi_{1,2}^{\pm}$, it follows that
$\int_0^{K} W_{YY}^2 d \eta \leq 1$. Lemma~\ref{lem:bounds} along
with
$$
W_X(X,\eta) =  \frac{1}{B} \zeta'\left(\frac{X}{B}\right) \left[\Phi_2(\eta) - 
\Phi_1(\eta)\right],
$$
implies that 
$$
\sup_{X,Y} |W_X(X,\eta)| \leq \frac{C}{B}.
$$
Combining this with the earlier estimate, we see that 
$$
\int_0^{M}  (V_X + U_Y + 2 W_X W_Y)^2 d \eta \leq C M \Delta_X^2 + \frac{C K^2}{B^2}
$$ 
Integrating this inequality in $X$, and using $M \geq K$, we
obtain the desired result.
\end{proof}
This Lemma provides a rigorous basis for our heuristic calculation for
the energies of the three regions at the beginning of this
section. Using this lemma, and the ideas form the heuristic
calculation, we have
\begin{theorem}
The minimizer of the elastic energy functional $\mathcal{I}$ in
(\ref{eq:unscaled}) subject to {\em compatible, no-stretch} boundary
conditions
$$
\begin{array}{lllc}
u = x, & v^{\pm} = \alpha^2 a
\Psi_{1,2}^{\pm}\left(\frac{y}{a_{1,2}}\right)+ y \pm \frac{1}{2}
\delta,& w^{\pm} = \sqrt{2} \alpha a \Phi_{1,2}^{\pm}\left(\frac{y}{a_{1,2}}\right),&
\mbox{ at } x = \pm L, \\ & & & \\ u = x & v^{\pm} = \pm L' \pm
\frac{\delta}{2}, & w^{\pm} = 0, & \mbox{ at } y = \pm L',
\end{array} 
$$
satisfies the upper bound
$$ \mathcal{I}^* \leq \min\left(C_1 \alpha^{7/3} \sigma^{5/3}
L^{1/3}, C_2 \alpha^2 \sigma^2 \frac{L}{a} \right)
$$ where $a = \min(a_1,a_2)$ and the constants $C_1$ and $C_2$ depend
only on $K$, where $K$ is the larger of the supports of
$\Phi^{\pm}_{1,2}$ divided by $a$. The constants are independent of
$a, \sigma, L,L'$.
\label{thm:non_slf_smlr}
\end{theorem}

\begin{proof}
We will consider two test solutions, and our upper bound will be the
minimum of the energies of the two test solutions.

One test solution is obtained by interpolating between the two
boundary conditions in the region $|y| \leq K a$, and
connecting this solution to the boundary conditions at $y = \pm L'$ by
a solution with a nearly uniform strain.

More precisely, we set $u = x$ and
$$ 
w^{\pm}(x,y) = \sqrt{2}\alpha a
\left[\zeta\left(\frac{1}{2}+\frac{x}{2L}\right)\Phi^{\pm}_2\left(\frac{y}{a}\right)
+ \left[1 - \zeta\left(\frac{1}{2}+\frac{x}{2
L}\right)\right]\Phi^{\pm}_1\left(\frac{y}{a}\right)\right]
$$ where $\zeta$ is as defined in
Lemma~\ref{lem:bndry_lyr}. Defining $\delta(x)$ by 
$$
\delta(x)= -\left[\frac{1}{2}\int_{-\infty}^0 \left[w^-_y(x,\eta)
\right]^2 d \eta +
\frac{1}{2}\int_0^{\infty} \left[w^+_y(x,\eta)\right]^2 d \eta + 2 \alpha w^+(x,0)\right]
$$
we set $v(x,y) = v_1(x,y) + v_2(x,y)$ where 
\begin{eqnarray}
v_1^{\pm}(x,y) & = &  \vartheta\left(\frac{y}{L'}\right)\left[-\frac{1}{2}\int_{\pm
\infty}^{y} \left[ w^{\pm}_y(x,\eta) \right]^2 d \eta + y \pm
\frac{1}{2}\delta(x)\right] \nonumber \\ 
v_2^{\pm}(x,y) & = & \left[1 - \vartheta\left(\frac{y}{L'}\right)\right] \left( y \pm \frac{1}{2}\delta\right) \nonumber
\end{eqnarray}
where $\vartheta$ is a smooth function such that $0 \leq \vartheta \leq 1$, 
$\vartheta(x) \equiv 1$ for $x \leq \frac{1}{3}$, and 
$\vartheta(x) \equiv 0$ for $x \geq \frac{2}{3}$.

We can estimate the energy of this test solution noting that
$\delta(x)$ is bounded by $C \alpha^4 a^2$ by
lemma~\ref{lem:bounds}, and using the appropriate rescaling of
lemma~\ref{lem:bndry_lyr} to bound the energy due to $v_1$. A
straightforward calculation gives 
$$
\mathcal{I} \leq  C \left[\frac{\alpha^4 a^5}{L^3} + \frac{\alpha^4 a^2 L'}{L} 
+ \frac{\sigma^2 \alpha^2 L}{a} + \frac{\sigma^2 \alpha^2 a^3}{L^3} +
\frac{\alpha^4 a^2 L}{L'} \right].
$$ This energy corresponds to $\mathcal{E}_I + \mathcal{E}_{III}$ in
our heuristic calculation, and this is of course reasonable, since we
do not have a region corresponding to the self-similar solution in
this test function. It is also clear that for $a \ll L \sim L'$, the
dominant term in the energy is $\sigma^2 \alpha^2 L/a$, and this will
be larger than the other terms provided that $\alpha^4 a^2 \ll
\sigma^2 \alpha^2 L/a$, {\em i.e} in the asymptotic regime $a \ll
\alpha^{-2/3} \sigma^{2/3} L ^{1/3}$.

Another test solution that we will consider is the following. Let
$\phi^{\pm},\psi^{\pm}$ be the zero asymptotic shift profile with
support $K_0$ that we constructed in Lemma~\ref{lem:chopped_off}.

We will first assume that $a_1 = a_2 = a$. Let $(\tilde{u},
\tilde{v}^{\pm}, \tilde{w}^{\pm})$ denote the self-similar solution
with profile $\phi^{\pm},\psi^{\pm}$, constructed in the same manner
as in Theorem.~\ref{thm:slf_smlr}, with $\rho(-1) = \rho(1) = a$.

We define the test solution by
\begin{eqnarray}
u(x,y) & = & \zeta\left(\frac{L-|x|}{b}\right) \tilde{u} + \left[ 1 -
\zeta\left(\frac{L-|x|}{b}\right)\right] x \nonumber \\ 
w^{\pm}(x,y) & = & \zeta\left(\frac{L - |x|}{b}\right) \tilde{w} + \sqrt{2}\left[1 -
\zeta\left(\frac{L-x}{B}\right)\right] \Phi^{\pm}_2\left(\frac{y}{a}\right) 
+ \left[1 - \zeta\left(\frac{L+x}{B}\right)\right] \Phi^{\pm}_1\left(\frac{y}{a}\right) 
\nonumber
\end{eqnarray}
where $\zeta$ is as defined in Lemma~\ref{lem:bndry_lyr}, and $b <
L$ is a length scale we will choose below. Defining $\delta(x)$ by
$$
\delta(x)= -\frac{1}{2}\left[\int_{-\infty}^0 \left[w^-_y(x,\eta)
\right]^2 d \eta +
\int_0^{\infty} \left[w^+_y(x,\eta)\right]^2 d \eta + 4 w^+(x,0)\right]
$$ it follows from the construction of $\tilde{w}$, and of $w$ that
$\delta(x) = 0$ if $L - |x| > b$. 

We set $v(x,y) = v_1(x,y) + v_2(x,y)$ where
\begin{eqnarray*}
v_1^{\pm}(x,y) & = & \vartheta\left(\frac{|y|}
{L'}\right)\left[-\int_{\pm \infty}^{y} \left[
w^{\pm}_y(x,\eta) \right]^2 d \eta + y \pm \frac{1}{2}\delta(x)\right]
\\ v_2^{\pm}(x,y) & = & \left[1 - \vartheta\left(\frac{|y|}{L'}\right)\right] \left( y \pm
\frac{1}{2}\delta\right)
\end{eqnarray*}
where $\vartheta$ is a smooth function such that $0 \leq \vartheta
\leq 1$, $\vartheta(x) \equiv 1$ for $x \leq \frac{1}{3}$, and
$\vartheta(x) \equiv 0$ for $x \geq \frac{2}{3}$. Note that $v(x,y) =
\tilde{v}(x,y)$ for $|y| \leq K_0 \alpha^{-1/3} \sigma^{1/3} L^{2/3}$
and $L - |x| \geq b$.

We can estimate the energy of this test solution as above and we obtain,
$$ 
\mathcal{I} \leq C(K) \left[\frac{\alpha^4 a^5}{b^3} +
\frac{\alpha^{11/3} a^2 \sigma^{1/3} L^{2/3}}{b} + \frac{\sigma^2
\alpha^2 b}{a} + \frac{\sigma^2 \alpha^2 a^3}{b^3} + \frac{\alpha^4
a^2 L}{L'} \right] + D \alpha^{7/3} \sigma^{5/3} L^{1/3}
$$ 
\begin{remark}
The constant $C$ depends on the boundary conditions through the
support $K$, but the constant $D$ is {\em independent} of $K$.
\end{remark}
This energy corresponds exactly to the energy $\mathcal{E}$ in our
heuristic calculation. Consequently, we get 
$$
\mathcal{I} \leq C\left(K,\frac{a}{\sigma}\right) \alpha^{7/3} \sigma^{5/3} L^{1/3}
$$ in the scaling regime $a/\sigma \sim O(1)$ for the choice $b =
\alpha^{5/6} \sigma^{2/3} L^{1/3}$.

$$ \mathcal{I} \leq D \alpha^{7/3} \sigma^{5/3} L^{1/3} +
C\left(K,\frac{a L^{2/3}}{\sigma^{5/3}}\right) \alpha^2 \sigma^2,
$$ in the scaling regime $a L^{2/3}\sigma^{-5/3} \sim O(1)$ for the
choice $b = a \sim \epsilon^{5/3} L ^{-2/3}$. $D$ is a {\em universal}
constant, independent of the scaling functions $\Phi$ and $\Psi$
determining the boundary conditions.

\end{proof}

\section{Structure of the minimal ridge} \label{sec:structure} 

In this section we will derive pointwise and integrated bounds for the
ridge width $\rho(X)$ and for the ridge sag $W(X,0)$. Lobkovsky's
analysis predicts that these two quantities should scale in the same
way, and further, the associated length scale is not uniform along the
ridge \cite{lobkovsky}.

Given a test solution $(u,v^{\pm},w^{\pm})$ that satisfies the
boundary conditions, we can naturally construct three different $x$
dependent length scales from the solution. They are, the inverse curvature 
$$
a(x) = \kappa^{-1}(x) = \left[ \int_{-\infty}^0 \left[w^-_{yy}(x,\eta)\right]^2 d \eta + \int^{\infty}_0 \left[w^+_{yy}(x,\eta)\right]^2 d \eta \right]^{-1}
$$
the ridge sag $w^+(x,0) = w^-(x,0)$, and the support
$$ k(x) = \sup\{y > 0 | [w^{\pm}(x,\pm y)]^2 + [v^{\pm}(x,\pm y)-y \mp
\delta/2]^2 + [u(x,\pm y)]^2 > 0\}
$$ 
In our construction for the self-similar test solutions yielding
the upper bound, we see that $a(x), w(x,0)$ and $k(x)$ all scale in
the same way as $\rho(X) = \rho(x/L)$. This suggests that the
structure of the ridge in the $y$-direction is given by a single
length scale, which depends on $x$, and further, this length scale is
given by our assumed scaling for $\rho(X)$ in the construction for the
upper bound. 

Our basic tool will be bounding the stretching energy using only the
length scales $a(x)$ and $w(x,0)$ at a given point $x$. Combining
these estimates with our estimates for the energy, we obtain
pointwise bounds for $a(x)$ and $w(x,0)$. This estimate for the
stretching energy is obtained in the following lemma.

\begin{lemma} Given the profile $W(X,.)$ at a point $X \in [-1,1]$, the 
stretching energy $E_s$ is bounded from below by
\begin{eqnarray}
\sqrt{E_s} & \geq & \max_{X \in [-1,1], Y \in \mathbb{R}, \delta \in
(0,1)} \, \frac{C (1 - \delta)}{1-X^2} \, \left[\sqrt{Y} \left(W(X,0)
+ \frac{Y}{\sqrt{2}}\right)^2 \right. \nonumber \\   
& & \left. + \frac{Y^{5/2}}{6} \left(1 - \frac{Y}{2 \delta \rho(X)}\right) - \frac{C'
A^3}{\sqrt{Y}}\right], \nonumber
\end{eqnarray}
where 
$$ 
\rho(X) = \left[\int_0^{\infty} {W^+_{YY}}^2(X,Y) dY + 
\int^0_{-\infty} {W^-_{YY}}^2(X,Y) dY\right]^{-1}.
$$
\label{lem:stretch}
\end{lemma}

\begin{proof}
The stretching energy is given by
$$
E_s = \int_{-\infty}^{\infty}\frac{1}{2} 
\left(\int_{-1}^1 W_X^2 dX \right)^2 dY 
$$
Since the integrand is non-negative, we have
\begin{eqnarray}
E_s & \geq & \max_Y  \int_{-Y}^{Y}\frac{1}{2} 
\left(\int_{-1}^1 W_X^2 dX \right)^2 dY \nonumber \\
& \geq & \max_Y \frac{1}{2Y} \left[\int_{-Y}^{Y}\int_{-1}^1 W_X^2 dX dY\right]^2, 
\label{ineq1}
\end{eqnarray}
by Jensen's inequality. For each $\eta \in [-Y,Y]$, and for all $X \in (-1,1)$, we have the
elementary inequality
\begin{eqnarray}
\int_{-1}^{1} W_X^2 (\xi,\eta) d \xi & \geq & \frac{(W(X,\eta) -
W(-1,\eta))^2} {1 + X} + \frac{(W(X,\eta) - W(1,\eta))^2}{1 - X},
\nonumber\\
& \geq & \frac{W^2(X,\eta)}{2} \left[\frac{1}{1+X} +
\frac{1}{1-X}\right] - 2 [ W^2(-1,\eta) + W^2(1,\eta)] \nonumber \\
& = & \frac{W^2(X,\eta)}{1 - X^2} - 2[W^2(-1,\eta) + W^2(1,\eta)] 
\label{ineq2}
\end{eqnarray}
$W(\pm 1,\eta)$ are given by the boundary conditions on the
frame. Using the scaling form for the boundary conditions, we see that
\begin{equation}
2 \left[\int_{-Y}^{Y} W^2(-1,\eta) d \eta + \int_{-Y}^{Y} W^2(1,\eta) d \eta \right]
\leq C A^3
\label{ineq3}
\end{equation}
All we now need is an estimate for $\iint W(\xi,\eta) d \xi d
\eta$. Let $\rho^{\pm}(X)$ be given by
$$
\rho^{\pm}(X) = \pm \left[\int_0^{\pm \infty} {W^{\pm}_{YY}}^2(X,Y) dY\right]^{-1}, 
$$ so that $[\rho(X)]^{-1} = [\rho^+(X)]^{-1} + [\rho^-(X)]^{-1}$. By
Lemma~\ref{lem:ineq}, we obtain
\begin{eqnarray} 
\int_0^{Y} [W^+(X,\eta)]^2 d \eta & \geq & \max_{Z \leq Y, \delta \in
(0,1)} \, (1 - \delta) \, \left[Z\left(W(X,0) + \frac{\beta^+(X)
Z}{2}\right)^2 \right. \nonumber \\
& & \left. + \frac{Z^3}{12} \left([\beta^+(X)]^2 - \frac{Z}{\delta
\rho^+(X)}\right)\right]. \nonumber 
\end{eqnarray}
Adding the corresponding result for $W^-$, we get
\begin{eqnarray} 
\int_{-Y}^{Y} [W(X,\eta)]^2 d \eta & \geq & \max_{Z_2 \leq Y, \delta \in
(0,1)} \, (1 - \delta) \, \left[Z_2\left(W(X,0) + \frac{\beta^+(X)
Z_2}{2}\right)^2 \right. \nonumber \\
& & \left. + \frac{Z_2^3}{12} \left([\beta^+(X)]^2 - \frac{Z_2}{\delta
\rho^+(X)}\right)\right] \nonumber \\
& & + \max_{Z_1 \leq Y, \delta \in
(0,1)} \, (1 - \delta) \, \left[Z_1\left(W(X,0) - \frac{\beta^-(X)
Z_1}{2}\right)^2 \right. \nonumber \\
& & \left. + \frac{Z_1^3}{12} \left([\beta^-(X)]^2 - \frac{Z_1}{\delta
\rho^+(X)}\right)\right] \nonumber 
\end{eqnarray}
By the matching condition, $\beta^+(X) = \beta^-(X) + \sqrt{2}$. This also
implies that $[\beta^-(X)]^2 + [\beta^+(X)]^2 \geq 1$. Replacing the
separate maximization over $Z_1 \leq Y$ and $Z_2 \leq Y$, by a single
maximization over $Z = Z_1 = Z_2 \leq Y$, we obtain
\begin{eqnarray} 
\int_{-Y}^{Y} [W(X,\eta)]^2 d \eta & \geq & \max_{Z \leq Y, \delta \in
(0,1)} \, 2 (1 - \delta) \, \left[Z\left(W(X,0) + \frac{Z}{\sqrt{2}}\right)^2 
+ \frac{Z^3}{12} \left(2 - \frac{Z}{\delta
\rho(X)}\right)\right] \nonumber
\end{eqnarray}
Combining this with inequalities (\ref{ineq1}),(\ref{ineq2}) and
(\ref{ineq3}) proves the lemma.
\end{proof}

We can now prove pointwise upper bounds for $\rho(X)$ and $W(X,0)$.
\begin{theorem} 
$(U,V^{\pm},W^{\pm})$ is a test solution that satisfies the boundary and 
the matching conditions. Also, $I(U,V^{\pm},W^{\pm}) \leq \bar{I}$. Then, 
$\exists  C_1,C_2,C_1',C_2'$ and $C_3'$ that only depend on $\bar{I}$ such 
that
$$
\rho(X) \leq C_1(1-X^2)^{2/5} + C_1'A.
$$ 
Also, the ridge sag satisfies the pointwise bound
$$ 
W^2(X,0) \leq \max\left(C_1' \rho^2(X), C_2' A^2 \left(\frac{A}{\rho(X)}\right)^{1/4}, 
C_3' \left(\frac{(1-X^2)^6}{\rho(X)}\right)^{1/7}\right).
$$ 
\label{thm:pointwise} 
\end{theorem}

\begin{proof}
Setting $\delta = \frac{1}{2}, Y = \frac{5}{7} \rho(X)$, Lemma~\ref{lem:stretch}
yields
$$ \sqrt{\bar{I}} \geq \sqrt{E_s} \geq \frac{C}{\sqrt{\rho(X)}(1-X^2)}
\left\{ [\rho(X)]^3 - C'A^3\right\}.
$$ If $\rho(X) \geq (2C')^{1/3} A$, it follows that $\rho^3(X) - C'A^3
\geq \frac{1}{2}\rho^3(X)$. Using this in the above inequality, we see
that either $\rho(X) < (2C')^{1/3} A$, or 
$$
\sqrt{\bar{I}} \geq \frac{C}{2 (1-X^2)} [\rho(X)]^{5/2}.
$$ Combining these two inequalities, we see that $\exists C_1,C_1'$
{\em only depending} on $\bar{I}$, such that 
$$
\rho(X) \leq C_1(1-X^2)^{2/5} + C_1'A.
$$ 
\end{proof}
This finishes the first part of the proof. To illustrate the pointwise
bounds for $W^2(X,0)$, we begin with a heuristic calculation. If
$W^2(X,0) \ll \rho^2(X)$, the dominant balance in
Lemma~\ref{lem:stretch} is
$$
Y^{5/2} \sim \frac{Y^{7/2}}{\rho}
$$ and this gives the characteristic scale $\tilde{Y} \sim \rho$. This
is the same calculation as above, and this gives $W^2 \ll \rho^2 \leq
(1-X^2)^{2/5}$ as in the previous part. If $W^2(X,0) \gg \rho^2(X)$,
after ignoring the boundary term $C'A^3 Y^{-1/2}$, the four remaining
terms in Lemma~\ref{lem:stretch} are of orders
$$ 
W^2(X,0) \sqrt{\tilde{Y}}, \quad W(X,0) \tilde{Y}^{3/2}, \quad \tilde{Y}^{5/2},
\quad \mbox{and} \quad \frac{\tilde{Y}^{7/2}}{\rho(X)},
$$ 
respectively. The dominant balance is between the first and the
last terms, and this gives the characteristic scale $\tilde{Y} \sim
[W^2(X,0) \rho(X)]^{1/3}$. Also, this gives 
$$ 
\sqrt{E_s} \sim \frac{C}{1-X^2} W^2(X,0) \sqrt{\tilde{Y}} \sim
\frac{C}{1-X^2} W^{7/3}(X,0) \rho^{1/6}(X).
$$
Rearranging, we get 
$$
W^2(X,0) \lesssim \left(\frac{(1-X^2)^6}{\rho(X)}\right)^{1/7}.
$$

We will now make these considerations precise.

\begin{lemma} 
$\exists C < \infty, C' > 0$ such that $|W(X,0)| \geq C \rho(X)$ implies that 
$$ 
\sqrt{Y} \left(W(X,0) + \frac{Y}{2}\right)^2 + \frac{Y^{5/2}}{6}
\left(1 - \frac{Y}{\rho(X)}\right) \geq C' W^{7/3}(X,0) \rho^{1/3}(X),
$$
where 
$$
Y = \left(\frac{6}{7}\right)^{1/3} \left[ W^2(X,0) \rho(X)\right]^{1/3}.
$$
\label{lem:easy}
\end{lemma}

\begin{proof} Let $W(X,0) = C_1 \rho(X)$, and $Y$ be as defined above. 
A direct calculation shows that
\begin{eqnarray}
\sqrt{Y} \left(W(X,0) + \frac{Y}{2}\right)^2 + \frac{Y^{5/2}}{6}
\left(1 - \frac{Y}{\rho(X)}\right) & \geq & 
\sqrt{Y} \left[ W^2 - |W| Y - \frac{Y^3}{6 \rho(X)}\right] \nonumber \\ 
& = & \left[\left(\frac{6}{7}\right)^{7/6} - C_1^{-1/3}\right] W^{7/3} \rho^{1/3} \nonumber
\end{eqnarray}
By taking $C_1$ sufficiently large ($> C$), we obtain the required
inequality. In particular, we can take $C = 8$ and $C' = 1/3$.
\end{proof}

We can now finish the proof of Theorem~\ref{thm:pointwise}.
\begin{proof}
If $W^2(X,0) \geq C^2 \rho(X,0)$, combining Lemma~\ref{lem:stretch}
and Lemma~\ref{lem:easy}, we obtain
$$
\sqrt{E_s} \geq \frac{C}{1-X^2} \left[\left(W^{14}(X,0) \rho(X)\right)^{1/6} 
- \frac{C' A^3}{\left(W^{2}(X,0) \rho(X)\right)^{1/6}}\right].
$$ 
If 
$$
\left(W^{14}(X,0) \rho(X)\right)^{1/6} \geq \frac{2 C' A^3}{\left(W^{2}(X,0) 
\rho(X)\right)^{1/6}},
$$
it follows that 
$$
\sqrt{\bar{I}} \geq \sqrt{E_s} \geq \frac{C}{2(1-X^2)} \left(W^{14}(X,0) \rho(X)\right)^{1/6}.
$$ 
Combining all the above considerations, we see that one of the
following inequalities has to hold
\begin{eqnarray}
W^2(X,0) & \leq & C_1' \rho^2(X), \nonumber \\
W^2(X,0) & \leq &  C_2' A^2 \left(\frac{A}{\rho(X)}\right)^{1/4}, \nonumber \\ 
W^2(X,0)& \leq & C_3' \left(\frac{(1-X^2)^6}{\rho(X)}\right)^{1/7}. \nonumber
\end{eqnarray}
This proves the theorem.
\end{proof}

\begin{remark} Our pointwise bounds are not optimal for $\rho(X)$. In particular, 
our construction for the upper bound shows that we can have test solutions with  
$$
\rho(X) \leq  C \max\left(A, \epsilon^{-1/3} (1- |X|), (1- X^2)^{2/3}\right),
$$ that have a uniformly bounded energy. 

For $A \ll (1 - |X|) \ll 1$, this function is asymptotically (in
$\epsilon$) much smaller than the upper bound from
Theorem~\ref{thm:pointwise}. In the remaining range, {\em i.e}, for
$1-|X| \lesssim A$, of for $C \lesssim 1-|X|$, where $C$ is an $O(1)$
constant, our pointwise upper bound captures the behavior of $\rho(X)$
in the self-similar construction.
\end{remark} 

\begin{remark} If $\rho(X)$ did scale like the pointwise upper bound in 
Theorem~\ref{thm:pointwise}, {\em i.e} $\rho(X) \sim C_1 A + C_2
(1-X^2)^{2/5}$, then the pointwise bounds for the ridge sag imply that 
$$
|W(X,0)| \leq C_1' A+ C_2'(1-X^2)^{2/5} \sim \rho(X).
$$
\end{remark}

We cannot obtain lower bounds for $\rho(X)$ or $W(X,0)$ purely by
energetic arguments as the following ``pinching'' argument shows. For
any given point $x \in (-L,L)$, we can consider the test solution
obtained by pinching at $x$, {\em i.e}, we set 
$$ u(x,s) = 0, \quad w(x,s) = \sqrt{2} \alpha a
\Phi\left(\frac{s}{a}\right), \quad v(x,s) = \alpha^2 a
\Psi\left(\frac{s}{a}\right) + s,
$$ where $a$ is a length scale that we are free to choose subject to
$a > \sigma \exp(- \epsilon^{-1/3})$ and $(\Phi,\Psi)$ is a no-stretch
profile with zero asymptotic shift that we constructed in
lemma~\ref{lem:chopped_off}. We use our construction for the
upper bound to construct minimal ridge solutions in $[-L,x]$ and
$[x,L]$.

These solutions connect smoothly at $x$, since our constructions have
$u(x',s) = u(x,s), v(x',s) = v(x,s)$ and $w(x',s) = w(x,s)$ for
sufficiently small $|x'-x|$. Also, $(L+x)^p + (L - x)^p \leq
2^{p+1}L^p$, for all $p > 0$. Combining this with our upper bounds, we
see that the energy of the ``pinched'' solution scales in the same
manner as the upper bound. Since the length scale $a$ can be chosen
(essentially) arbitrarily small, it follows that the energetics are
not enough to give us a pointwise lower bound on the ridge-sag
$W(X,0)$ or the the ridge-width $\rho(X)$.

We need these pointwise lower bounds to prove rigorous scaling results
for $\rho(X)$ and $W(X,0)$.To obtain such results, we need to use the
fact that the solution of interest is a minimizer, {\em i.e}, the
first variations vanish. This type of analysis is carried out for the
structure of an Austenite-Martensite boundary in \cite{KM,conti}. A
similar analysis is possible for the minimal ridge, and we will
present the details elsewhere.

\section{Discussion} \label{sec:discussion}

We conclude our discussion by indicating some of the issues/open
problems relating to thin elastic sheets in general and to minimal
ridges in particular, and in the process we point out the relevance of
our results to some of these questions.

We have proved rigorous scaling laws for the energy of a single
minimal ridge with a geometrically linear FvK {\em ansatz}. A natural
question is the extension of these results to the Nonlinear FvK
energy, and also to full three dimensional elasticity, as in
\cite{BCDM02}. It is easy to extend this to a mixed energy functional
where the bending energy is treated in a geometrically linear fashion,
but we use the ``full'' energy $W_{2D}$ for the in-plane
stretching. Extending this analysis to the Geometrically nonlinear
functional, or to full three dimensional elasticity will require some
new techniques \cite{BCDM02}.

Another problem is to show that the scaling laws also hold for a real
crumpled sheet, where the forcing is not through clamping the
boundaries to a frame, but rather through the confinement in a small
volume. In this case, there are interesting global geometric and
topological considerations, some of which are explored in
Refs.~\cite{immersion_thm} and \cite{high_d}. As Lobkovsky and Witten
\cite{LobWit} argue, the boundary condition that the deformation goes
to zero far away from the ridge implies that the ridges do not
interact with each other significantly. The ridges can be considered
the {\em elementary excitations} of a crumpled sheet.

More precisely, we have constructed ridge solutions with zero
asymptotic shift, that are exactly strain free on the
boundaries. These solutions give ``non-interacting'' ridges and
patching these solutions together, it is possible to construct test
solutions for a sheet confined inside a sphere. This gives us upper
bounds which scale in the same way as the energy of a single ridge.

To show the corresponding lower bound, we have to show that
confinement actually leads to the formation of ridges, and that the
competition between the bending and the stretching energy for this
situation has the same form as in lemma~\ref{lem:essential}. In this
context, we expect that global topological considerations, as well as
the non self-intersection of the sheet will play a key role in the
analysis, as they do in the analysis of elastic rods (one dimensional
objects) \cite{GMSvdM}.

As we note above, the blistering problem is described by the same
elastic energy (Eq.~(\ref{eq:unscaled})), but with different boundary
conditions. Our results show an interesting contrast with results for
the blistering problem. Ben Belgacem {\em et al.} have shown that
\cite{BCDM}, for an isotropically compressed thin film, the energy of
the minimizer satisfies
$$
c \lambda^{3/2} \sigma L \leq {\mathcal I}
\leq C \lambda^{3/2} \sigma L,
$$
where $L$ is a typical length scale of the domain, and $\lambda$ is
the compression factor. A construction for the upper bound strongly
suggests that the minimizers develop an infinitely branched network
with oscillations on increasingly finer scales as $\sigma
\rightarrow 0$. In contrast, our results indicate that the energy of a 
minimal ridge satisfies
$$
c \sigma^{5/3} L^{1/3} \leq {\mathcal I}
\leq C \sigma^{5/3} L^{1/3},
$$ and the energy concentrates in a region of width
$\sigma^{1/3}L^{2/3}$. This shows that the nature of the solution of
the variational problem for the elastic energy in (\ref{eq:unscaled})
depends very strongly on the boundary conditions. In particular the
very nature of the energy minimizers is different for the two problems
-- For the blistering problems, as $\sigma \rightarrow 0$ the
minimizers develop a branched network of folds refining towards the
boundary.

Finally,one would like to prove $\Gamma$--convergence and find the
$\Gamma$--limit \cite{DalMaso,DeGDM} for the elastic energy as $\sigma
\rightarrow 0$. The difference in the scaling of the energy minimum
for the minimal ridge, and the blistering problem shows that the
$\Gamma$--limit of the elastic energy depends crucially on the imposed
boundary conditions.

The analysis in this paper only pertains to situations where the
configuration of the sheet is either smooth, or consists of a finite
number of minimal ridges.  More precisely, the sheet configurations
$\phi:\mathcal{S} \rightarrow \mathbb{R}^3$ is piecewise smooth,
strain-free a.e., has gradient $D \phi$ in $BV$, and the singular
support of $D^2 \phi$ lives on a finite union of straight line
segments. For the boundary conditions that admit such configurations,
our analysis suggests that the asymptotic energy is on the scale
$\sigma^{5/3}$ as $\sigma \rightarrow 0$. Further, if the following
limit exists, the is necessarily follows that
$$ 
\bar{\mathcal{I}}[\phi] = \Gamma-\lim_{\sigma \rightarrow 0}
\sigma^{-5/3} \mathcal{I} = C \sum \alpha_j^{7/3} l_j^{1/3},
$$ where $l_j$ is the length of the $j$th segment in the singular
support of $D^2 \phi$, and $\alpha_j$ is the jump in $D \phi$ across
the segment $\alpha_j$. Note that, because of the $l_j^{1/3}$
dependence, the $\Gamma$--limit {\em cannot be written} as the
integral with respect to the $\mathcal{H}^1$ measure on the singular
support of $D^2 \phi$, of a {\em local energy density}, which only
depends on $\alpha_j$.

The reason we get a $l_j^{1/3}$ dependence instead of a linear
dependence in $l_j$ is that the ridges have a nonuniform structure
along their length for any $\sigma > 0$. Any formulation of the
$\Gamma$--limit should therefore incorporate a ``hidden'' variable
which reflects the nonuniform structure, although, this non-uniformity
is no longer detectable in the limiting $\sigma \rightarrow 0$
configuration. A natural candidate for this variable is the scaled
ridge width $\rho(x)$ (= the inverse curvature $a(x)$ defined as in
Sec.~\ref{sec:structure}), where $x$ is a coordinate along a ridge. We
expect that the $\Gamma$--limit can be written as an integral of an
energy density with respect to $\mathcal{H}^1$ measure on the defect
set, where the energy density depends on $\alpha$ and also on
$\rho(x)$ and derivatives $\rho'(x)$ and $\rho''(x)$. In fact,
Eq.~(\ref{eq:estimates}) strongly suggest that the $\Gamma$-limit for
the energy of a single ridge can be written as
$$ \bar{\mathcal{I}}[\phi] = C \sum \alpha_j^{7/3}
\inf_{\rho(.)}\int_0^{l_j} \left[\left([\rho''(x) \rho(x)]^2 +
[\rho'(x)]^4\right) \rho(x) + \frac{1}{\rho(x)} \right]dx
$$ where $x$ is a coordinate along the ridge, and the infimum is over
all smooth functions $\rho$ that vanish at both the endpoints of the
ridge.

We hope this paper spurs further investigation of the question of the
$\Gamma$--limit of the elastic energy. This question is very much
open, there are no proofs for either $\Gamma$--convergence, or of our
conjectured structure for the $\Gamma$--limit.

\section*{Acknowledgements}

I am very grateful to Brian DiDonna, Felix Otto, Stefan M\"{u}ller and
Tom Witten for many helpful and enlightening conversations on the
subject of crumpled sheets. This work was supported in part by the
National Science Foundation through its MRSEC program under Award
Number DMR-0213745, and by a NSF CAREER Award DMS-0135078. This work
is also supported in part by a Research Fellowship from the Alfred
P. Sloan Jr. Foundation.

\providecommand{\bysame}{\leavevmode\hbox to3em{\hrulefill}\thinspace}


\begin{thebibliography}{10}

\bibitem{Audoly}
B.~Audoly, \emph{Stability of straight delamination blisters}, Phys. Rev. Lett.
  \textbf{83} (1999), 4124.

\bibitem{pomeau}
M.~Ben~Amar and Y.~Pomeau, \emph{Crumpled paper}, Proc. Roy. Soc. London Ser. A
  \textbf{453} (1997), 729.

\bibitem{BCDM}
H.~Ben~Belgacem, S.~Conti, A.~DeSimone, and S.~M{\"u}ller, \emph{Rigorous
  bounds for the {F}\"oppl-von {K}\'arm\'an theory of isotropically compressed
  plates}, J. Nonlinear Sci. \textbf{10} (2000), no.~6, 661--683.

\bibitem{BCDM02}
\bysame, \emph{Energy scaling of compressed elastic films - three-dimensional
  elasticity and reduced theories}, Arch. Rat. Mech. Anal. \textbf{164} (2002),
  1--37.

\bibitem{BPCB00}
A~Boudaoud, P~Patricio, Y~Couder, and M~Ben~Amar, \emph{Dynamics of
  singularities in a constrained elastic plate}, Nature \textbf{407} (2000),
  no.~6805, 718--720.

\bibitem{vertex}
E.~Cerda, S.~Cha\"{\i}eb, F.~Melo, and L.~Mahadevan, \emph{Conical dislocations
  in crumpling}, Nature \textbf{401} (1999), 46.

\bibitem{maha}
E.~Cerda and L.~Mahadevan, \emph{Concical surfaces and crescent singularities
  in crumpled sheets}, Phys. Rev. Lett. \textbf{80} (1998), 2354.

\bibitem{dcone_exp}
S.~Cha\"{\i}eb, F.~Melo, and J-C. G\'eminard, \emph{Experimental study of
  developable cones}, Phys. Rev. Lett. \textbf{80} (1998), 2358.

\bibitem{ciarlet2}
Philippe~G. Ciarlet, \emph{Mathematical elasticity. {V}ol. {I}{I}},
  North-Holland Publishing Co., Amsterdam, 1997, Theory of plates.

\bibitem{conti}
Sergio Conti, \emph{Branched microstructures: scaling and asymptotic
  self-similarity}, Comm. Pure Appl. Math. \textbf{53} (2000), no.~11,
  1448--1474.

\bibitem{dacorogna}
Bernard Dacorogna, \emph{Direct methods in the calculus of variations},
  Springer-Verlag, Berlin, 1989.

\bibitem{DalMaso}
Gianni Dal~Maso, \emph{An introduction to {$\Gamma$}-convergence}, Birkh\"auser
  Boston Inc., Boston, MA, 1993.

\bibitem{DeGDM}
Ennio De~Giorgi and Gianni Dal~Maso, \emph{${\Gamma}$-convergence and calculus
  of variations}, Mathematical theories of optimization (Genova, 1981),
  Springer, Berlin, 1983, pp.~121--143.

\bibitem{ridge_buckling}
B.~A. DiDonna, \emph{Scaling of the buckling transition of ridges in thin
  sheets}, Phys. Rev. E \textbf{66} (2002), no.~1, 016601.

\bibitem{ridge_strength}
B.~A. DiDonna and Witten~T. A., \emph{Anomalous strength of membranes with
  elastic ridges}, Phys. Rev. Lett. \textbf{87} (2001), no.~20, 206105.

\bibitem{high_d}
B.~A. DiDonna, T.~A. Witten, S.~C. Venkataramani, and E.~M. Kramer,
  \emph{Singularities, structures, and scaling in deformed $m$-dimensional
  elastic manifolds}, Phys. Rev. E (3) \textbf{65} (2002), no.~1, part 2,
  016603, 25.

\bibitem{evans}
Lawrence~C. Evans, \emph{Weak convergence methods for nonlinear partial
  differential equations}, Published for the Conference Board of the
  Mathematical Sciences, Washington, DC, 1990.

\bibitem{evans2}
\bysame, \emph{Partial differential equations}, American Mathematical Society,
  Providence, RI, 1998.

\bibitem{GDeSOC02}
G.~Gioia, A.~DeSimone, M.~Ortiz, and A.~M. Cuiti\~{n}o, \emph{Folding
  energetics in thin-film diaphragms}, R. Soc. Lond. Proc. Ser. A Math. Phys.
  Eng. Sci. \textbf{458} (2002), no.~2021, 1223--1229.

\bibitem{GMSvdM}
O.~Gonzalez, J.~H. Maddocks, F.~Schuricht, and H.~von~der Mosel, \emph{Global
  curvature and self-contact of nonlinearly elastic curves and rods}, Calc.
  Var. Partial Differential Equations \textbf{14} (2002), no.~1, 29--68.

\bibitem{KJ}
W.~Jin and R.~V. Kohn, \emph{Singular perturbation and the energy of folds}, J.
  Nonlinear Sci. \textbf{10} (2000), no.~3, 355--390.

\bibitem{JS}
Weimin Jin and Peter Sternberg, \emph{Energy estimates for the von {K}\'arm\'an
  model of thin-film blistering}, J. Math. Phys. \textbf{42} (2001), no.~1,
  192--199.

\bibitem{JS1}
\bysame, \emph{In-plane displacements in thin film blistering}, preprint, 2001.

\bibitem{KM}
Robert~V. Kohn and Stefan M{\"u}ller, \emph{Surface energy and microstructure
  in coherent phase transitions}, Comm. Pure Appl. Math. \textbf{47} (1994),
  no.~4, 405--435.

\bibitem{eric}
E.~M. Kramer and T.~A. Witten, \emph{Stress condensation in crushed elastic
  manifolds}, Phys. Rev. Lett. \textbf{78} (1997), 1303.

\bibitem{landau_elastic}
L.~D. Landau and E.~M. Lifshitz, \emph{Theory of elasticity}, Pergamon Press,
  London, 1959.

\bibitem{science.paper}
A.~Lobkovsky, S.~Gentges, H.~Li, D.~Morse, and T.~A. Witten, \emph{Scaling
  properties of stretching ridges in a crumpled elastic sheet}, Science
  \textbf{270} (1995), 1482.

\bibitem{lobkovsky}
A.~E. Lobkovsky, \emph{Boundary layer analysis of the ridge singularity in a
  thin plate}, Phys. Rev. E. \textbf{53} (1996), 3750.

\bibitem{LobWit}
A.~E. Lobkovsky and T.~A. Witten, \emph{Properties of ridges in elastic
  membranes}, Phys. Rev. E \textbf{55} (1997), 1577.

\bibitem{love}
A.~E.~H. Love, \emph{A treatise on the {M}athematical {T}heory of
  {E}lasticity}, Dover Publications, New York, 1944, Fourth Ed.

\bibitem{MB02}
T.~Mora and A.~Boudaoud, \emph{Thin elastic plates: On the core of developable
  cones}, Europhys. Lett. \textbf{59} (2002), no.~1, 41--47.

\bibitem{GiOrtz}
Michael Ortiz and Gustavo Gioia, \emph{The morphology and folding patterns of
  buckling-driven thin-film blisters}, J. Mech. Phys. Solids \textbf{42}
  (1994), no.~3, 531--559.

\bibitem{eran}
E.~Sharon, B.~Roman, M.~Marder, G-S. Shin, and H.~L. Swinney, \emph{Buckling
  cascade in free thin sheets}, Nature \textbf{419} (2002), 579.

\bibitem{eran2}
E.~Sharon, S.~Smith, H.~L. Swinney, and B.~Roman and, \emph{Geometrically
  induced buckling cascade in free sheets}, preprint, 2003.

\bibitem{immersion_thm}
S.~C. Venkataramani, T.~A. Witten, E.~M. Kramer, and R.~P. Geroch,
  \emph{Limitations on the smooth confinement of an unstretchable manifold}, J.
  Math. Phys. \textbf{41} (2000), no.~7, 5107--5128.

\bibitem{ridge1}
Shankar~C. Venkataramani, \emph{Lower bounds for the energy in a crumpled
  elastic sheet -- a minimal ridge}, Submitted to Nonlinearity, Preprint {\tt
  http://arxiv.org/abs/math.AP/0210012}.

\bibitem{Young}
L.~C. Young, \emph{Lectures on the calculus of variations and optimal control
  theory}, W. B. Saunders Co., Philadelphia, 1969.

\bibitem{ziemer}
William~P. Ziemer, \emph{Weakly differentiable functions}, Graduate Texts in
  Mathematics, vol. 120, Springer-Verlag, New York, 1989, Sobolev spaces and
  functions of bounded variation.

\end{thebibliography}
\end{document}